%% file: Scov_MROLiveira_arxiv_2021_v3.tex
\RequirePackage{fix-cm}
\documentclass{svjour3RO}                     % onecolumn (standard format)
\smartqed  % flush right qed marks, e.g. at end of proof
\usepackage{graphicx}
\usepackage{lipsum}
\usepackage{overpic}
\usepackage{subfig}
\usepackage{eqnarray}
\usepackage{float}
\usepackage{colortbl}

\usepackage{rotating}
\usepackage{pifont} 
\usepackage{commath}
\usepackage{xcolor}
\usepackage{multicol}
\usepackage{enumitem}
\setlist{itemsep=0.2pt}
\usepackage{textcomp} 
\usepackage{comment}
\usepackage{amsmath}  
\usepackage{graphicx}
\usepackage{booktabs,multirow}
\usepackage{caption}
\usepackage{fancyvrb}
\usepackage{amsfonts}
\usepackage{subfloat}
\usepackage{wrapfig}
\usepackage{upgreek}
\usepackage[square,sort,numbers]{natbib}
\usepackage{footnote}
\usepackage{pifont}%%%%Acrescentar no main file - Ro, 07/set/2018
\newcommand{\eqnref}[1]{(\ref{#1})}

\newcommand{\dps} {\displaystyle}

\renewcommand{\Re}{\mathbb{R}}

\newcommand{\boma}[1]{\mbox{\boldmath ${#1}$}}
\newcommand{\qedwhite}{\hfill \ensuremath{\Box}} 

 \journalname{}
\begin{document}

\title{Association measures for interval variables\thanks{This research has been supported by Funda\c{c}\~{a}o para a Ci\^{e}ncia
		e Tecnologia (FCT), Portugal, through the projects %UID/Multi/04621/2013,
		UIDB/04621/2020 and UIDP/04621/2020 of 286 CEMAT/IST-ID,
		UID/EEA/50008/2013,  
		PTDC/EGE-ECO/30535/2017, and 
		PTDC/EEI-TEL/32454/2017.}
%Grants or other notes
%%%%Measuring association with interval-valued data "
%%%%"Association measures for interval variables"
%------------------------------------------------------------
}


%\author{}
\author{M. Ros\'{a}rio Oliveira\and
	    Margarida Azeitona  \and
	    Ant\'{o}nio Pacheco  \and
	    Rui  Valadas %etc.
}

%-------------

%---------------

\authorrunning{M.~Ros\'{a}rio Oliveira \textit{et al.}} % if too long for running head

%\institute{}
\institute{M.~Ros\'{a}rio Oliveira, M. Azeitona, A. Pacheco \at
              CEMAT and Mathematics Department, Instituto Superior T\'{e}cnico, Universidade de Lisboa, Portugal \\
              \email{rosario.oliveira@tecnico.ulisboa.pt}           %  \\
           \and
           R.  Valadas \at
              Instituto de Telecomunica\c{c}\~{o}es and Dep. Electrical and Computer Engineering, Instituto Superior T\'{e}cnico, Universidade de Lisboa, Portugal
}

\date{January 2021}

\maketitle

\begin{abstract}
	Symbolic Data Analysis (SDA) is a relatively new field of statistics that extends conventional data analysis by taking into account intrinsic data variability and structure. Unlike conventional data analysis, in SDA the features characterizing the data can be multi-valued, such as intervals or histograms. SDA has been mainly approached from a sampling perspective. In this work, we propose a model that links the micro-data and macro-data of interval-valued symbolic variables, which takes a populational perspective. Using this model, we derive the micro-data assumptions underlying the various definitions of symbolic covariance matrices proposed in the literature, and show that these assumptions can be too restrictive, raising applicability concerns. We analyze the various definitions using worked examples and four datasets. Our results show that the existence/absence of correlations in the macro-data may not be correctly captured by the definitions of symbolic covariance matrices and that, in real data, there can be a strong divergence between these definitions. Thus, in order to select the most appropriate definition, one must have some knowledge about the micro-data structure.
%%%%%OLD version	
%Symbolic Data Analysis (SDA) is a relatively new field of statistics that extends conventional data analysis by taking into account intrinsic data variability and structure. Unlike conventional data analysis, in SDA the features characterizing the data can be multi-valued, such as intervals or histograms. SDA has been mainly approached from a sampling perspective and, in this work, we introduce population formulations of the symbolic mean, variance, covariance, correlation, covariance matrix, and correlation matrix for interval-valued symbolic variables, providing a theoretical framework that gives support to interval-valued SDA. Moreover, we provide an interpretation of the various definitions of covariance and correlation matrices according to the structure of micro-data, which allows selecting the model that best suits specific datasets. Our results are illustrated using two datasets. Specifically, we select the most appropriate model for each dataset using goodness-of-fit tests and quantile-quantile plots, and provide an explanation of the micro-data based on the covariance matrix.	
%%%%OLD version	

%
\end{abstract}

%-------------------------------------------------------
\input{Rev_Intro_SCov_v2}		  %%%Section 1.
%-------------------------------------------------------

%-------------------------------------------------------
\input{Rev_Reasoning_SCov_v1_APa_v4}   %%%Section 2.%%% - Change the name
%-------------------------------------------------------
 
 %-------------------------------------------------------
\input{Rev_State_of_the_art_v2} %%%%%Subsection -  Sec 2
%-------------------------------------------------------
%%%-------------------------------------------------------
\input{Rev_Null_SCov_v2}
%%%%%Subsection - CONCLUIDO
%%%-------------------------------------------------------

%-------------------------------------------------------
\input{Rev_Example_SCov_v2}	   %%%Section 4.%%% CONCLUIDO%%%%
%%-------------------------------------------------------

%-------------------------------------------------------
\input{Rev_Conclus_SCov_v3}	   %%%Section 5. %%% CONCLUIDO%%%%
%-------------------------------------------------------

%\section*{Acknowledgements}
%%This research has been supported by Funda\c{c}\~{a}o para a Ci\^{e}ncia
%%e Tecnologia (FCT), Portugal, through the projects UID/Multi/04621/2013, UID/EEA/50008/2013, and PTDC/EEI-TEL/5708/2014. 
%We thank the reviewers for their constructive comments which greatly enrich the paper.

\bibliographystyle{spmpsci}      % mathematics and physical sciences
\bibliography{Scov_MROLiveira_arxiv_2021_v3}

\end{document}

%% file: Rev_Intro_SCov_v2.tex
\section{Introduction}\label{sec:Introd}

The volume and complexity of available data in virtually all sectors of society has grown enormously, boosted by the globalization and the massive use of the Internet. New statistical methods are required to handle this new reality, and Symbolic Data Analysis (SDA), proposed by Diday in \cite{Diday1987}, is a promising research area.%RV 30/10/2020

When trying to characterize datasets, it may not be convenient to deal with the individual data observations (e.g. because the sample size is too large), or we may not have access to the individual data observations (e.g. because of privacy restrictions). In conventional data analysis, this problem is usually handled by providing single-valued summary statistics of the data characteristics (e.g. mean, variance, quantiles). The analysis can consider multiple characteristics, but these characteristics can only be single-valued. SDA extends conventional data analysis by allowing the description of datasets through multi-valued features, such as intervals, histograms, or even distributions \cite{BillardDiday2003,BillardDiday2006,Brito:2014}. These features are called \textit{symbolic variables}.%RV 30/10/2020

Suppose we want to analyze textile sector per country, e.g. in terms of two characteristics: number of customers and profit. Suppose also that we only have access to summary information per country, and not to the data of each individual company. Conventional data analysis can only deal with single-valued features, like the profit variance, profit mean, or the mean number of customers. Instead, in SDA, the features (the symbolic variables) can be multi-valued, e.g. one feature can be the minimum and maximum profits, and another can be a histogram of the number of customers.%RV 30/10/2020

One of the main benefits of SDA has to do with the way individual data characteristics (e.g. profit or number of customers) are described. In conventional data analysis, since only single-valued features are available, one may need many features to describe a given characteristic. Moreover, the features are treated in the same way, irrespective of the characteristic they represent. For example, one may create as features to characterize the profit of the textile sector the mean, the variance, the maximum, the minimum, the median, the first and third quartiles, and so on. There is then an inflation of features to explain a single data characteristic (the profit, in this case). SDA allows explaining single data characteristics through single symbolic data variables, better tailored to analyze that specific characteristic, and with potential gains in terms of dimensionality.%RV 30/10/2020

In SDA, the original data is called \textit{micro-data} and the aggregated data is called \textit{macro-data}. In the previous example, the micro-data would be the data of individual companies (labeled with the country they belong to), and the macro-data the interval of profit (between maximum and minimum) or the histogram of the number of customers, of the companies of each country. Our main interest in this paper is on interval-valued data \cite{Noirhomme.Brito:2011,Zhang.Sisson:2020}, where macro-data corresponds to the interval between minimum and maximum of micro-data values.%RV 30/10/2020

SDA is a relatively new field of statistics and has been mainly approached from a sampling perspective. The works  \cite{Bertrand2000,BillardDiday2006,Billard2008} introduced measures of location, dispersion, and association between symbolic random variables, formalized as a function of the observed macro-data values. The sample covariance (correlation) matrices were addressed in the context of symbolic principal component analysis in \citep{Chouakria1998,Le-Rademacher2008,Le-RademacherBillard2012,Wang2012,Vilela2015,Oliveira.et.al:2017} 
and more recently in factor analysis \cite{Cheira.et.al:2017}. In \cite{Oliveira.et.al:2017} the authors established relationships between several proposed methods of symbolic principal component analysis and available definitions of sample symbolic variance and covariance. Other areas of statistics have also been addressed by SDA like clustering  (e.g.\ \cite{Carvalho.Lechevallier:2009,Sato.Ilic:2011}), discriminant analysis (see e.g.\ \cite{DSilva.Brito:2015,Queiroz.et.al:2018}), regression analysis  (see e.g.\ \cite{LNeto.et.al:2011,Dias.Brito:2017}), and time series (see e.g.\ \cite{Maia.et.al:2008,Teles.Brito:2015}).%RV 30/10/2020

Parametric approaches for interval-valued variables have also been considered \cite{LeRademacher.Billard:2011, LNeto.et.al:2011, Brito2012,DSilva.Brito:2015, Dias.Brito:2017, DSilva.et.al:2017, Cheira.et.al:2017}. Authors in \cite{LeRademacher.Billard:2011} derived maximum likelihood estimators for the mean and the variance of three types of symbolic random variables: interval-valued, histogram-valued, and triangular distribution-valued variables. In \cite{LNeto.et.al:2011}, authors have formulated interval-valued variables as bivariate random vectors in order to introduce a symbolic regression model based on the theory of generalized linear models. The works \cite{Brito2012,DSilva.Brito:2015,DSilva.et.al:2017} have followed a different approach. In their line of work, the centers and the logarithms of the ranges are collected in a random vector with a multivariate (skew-)normal distribution, which is used to derive methods for the analysis of variance \cite{Brito2012}, discriminant analysis \cite{DSilva.Brito:2015}, and  outlier detection \cite{DSilva.et.al:2017} of interval-valued variables. More recently, the need to specify a micro-data model for the underlying macro-data was addressed by \cite{Zhang.Sisson:2020,Beranger.et.al:2020}. The authors constructed maximum likelihood functions for symbolic data based on how the macro-data is derived from the micro-data.%RV 30/10/2020

Despite of previous work, the area of SDA is lacking theoretical support and our work is a step in this direction. Preferably, the statistical methods of SDA should be grounded on populational formulations, as in the case of conventional methods. A populational formulation allows a clear definition of the underlying statistical model and its properties, and the derivation of effective estimation methods.%RV 30/10/2020

In this paper, we use population formulations of the sample symbolic mean vector, sample symbolic covariance matrices and sample correlation matrices available in the literature. To provide an interpretation of each definition, we propose a model that links the structure of the micro-data and of the macro-data. The model assumes that the micro-data associated with a certain random interval are not observable (latent), and has a mean equal to the mean of the interval center. Using the model, we derive the assumptions on the micro-data structure subjacent to each definition of symbolic covariance matrices. We also show that the symbolic correlations proposed in the literature are quantities between -1 and 1, as in the conventional case, independently of the relationship between the micro- and macro-data. Using worked examples, we discuss the meaning of null covariance, and show cases where the existence/absence of correlation in the macro-data is not captured by the definitions of symbolic correlation matrices. Finally, we explore the various definitions using real data examples. When micro-data is available, we select the most appropriate definition for each dataset using goodness-of-fit tests and quantile-quantile plots, and provide an explanation of the micro-data based on the covariance matrices. However, the results show that there can be a large divergence between definitions, meaning that, in general, one needs some information on the micro-data structure to decide about the most appropriate definition.%RV 2/11/2020

The paper is structured as follows. In Section \ref{sec:Reasoning} we introduce the model linking the micro-data with the macro-data of interval-valued symbolic variables and discuss several worked examples. Section \ref{sec:Examples} introduces the real data examples. Finally, Section \ref{sec:conclusions} presents the conclusions of the paper and gives some directions for future work.%RV 2/11/2020

%% file: Rev_Reasoning_SCov_v1_APa_v4.tex
\section{Symbolic means, variances, and covariances}
\label{sec:Reasoning}

In this work, we focus on the study of interval-valued variables \cite{Diday2000} and  interval-valued random vectors, as next defined.
\begin{definition}  
\label{def:IntervalVariable}
	$X=[A,B]$ is an interval-valued  random variable defined on the probability space $(\Omega, {\cal F},P)$ if and only if $A$ and $B$ are random variables defined on $(\Omega, {\cal F},P)$ such that $P(A \leq B)=1$.
	A $p$-dimensional interval-valued random vector $\boma{X}$ is a vector  of $p$ interval-valued  random variables, $X_1$, $X_2$,\dots,$X_p$, all defined on the same probability space.
\end{definition}%RO 30/10/2020

We consider that an interval-valued random variable, $X=[A, B]$, besides being represented by the interval limits $A$ and $B$, is also represented by the center and range of the interval:
\begin{align}
C=\dps \frac{A + B}{2}\quad \text{and} \quad R=B - A.\label{eq:transCR}
\end{align}%RV 27/10/2017
Similarly a $p$-dimensional interval-valued random vector $\boma{X}=(X_1,X_2,\ldots,X_p)^t$ is characterized by the vector $( \boma{C}^t,\boma{R}^t)^t$  where $\boldsymbol C = (C_1, \dotsc, C_p)^t$ is the vector of centers and $\boldsymbol R = (R_1, \dotsc, R_p)^t$ the vector of ranges describing the object, i.e.,
\begin{align}
C_i=\dps \frac{A_i + B_i}{2}\quad \text{and} \quad R_i=B_i - A_i\label{eq:transCRi}
\end{align}%RV 27/10/2017
for $i=1,2,\ldots ,p$.

In the next subsection, we start by proposing a model that establishes a natural link between micro-data and macro-data. Afterwards, we derive the means, variances, and covariances of the micro-data; we refer to these quantities as symbolic means, symbolic variances, and symbolic covariances.

\subsection{A model linking micro-data with macro-data}
\label{subsec:Reasoning_MicroData_Models}

In this model, we define a random vector $\boma{A}$ representing micro-data, as a function of the random vectors of centers and ranges, $\boma{C}$ and $\boma{R}$, that characterize the associated macro-data, along with a random vector $\boma{U}$ that characterizes the structure of the micro-data given the associated macro-data. 

Note that a realization of $\boma{A}$ is a point in the hyper-rectangle associated with the random interval-valued vector $\boma{X}$, characterized by its center, $\boma{C}$, and range,  $\boma{R}$. 
In detail, we choose $\boma{A}=(A_{1},\ldots,A_{p})^t$ such that:
\begin{equation}
\label{Aij}
	A_{j}=C_j+U_{j}\frac{\dps R_j}{\dps 2},\; j=1,2,\ldots,p,
\end{equation}%RV 3/11/2017
\noindent  with the weights $U_{j}$ being random variables with support on the interval $[-1,1]$. 

We note that $U_{j} R_j/2$ [that is the term that is added to the center of the $j$-th component of the random macro-data, $C_j$, to produce the $j$-th component of an associated micro-data, namely $A_j$] is the deviation of the $j$-th component of the micro-data to the center of the interval of the $j$-th component of its associated macro-data $\left[C_j - R_j/2, C_j + \dps R_j/2\right]$.

By imposing conditions on  $\boma{U}$, we obtain the symbolic mean vector and the symbolic covariance matrix, i.e., the mean vector and the covariance matrix  of $\boma{A}$, denoted by $\boma{\mu}={\rm E}(\boma{A})$ and $\boma{\Sigma}={\rm Var}(\boma{A})$, respectively.  
In the remaining of this subsection  we consider the assumption next stated (Assumption 1), which in turn leads to Theorem~\ref{Teo:Assumption_Cov_matrices}.\\

\noindent  {\bf Assumption 1}: The random vector $\boma{U}= (U_{1},U_2,\ldots,U_{p})^t$ has zero mean and is independent of the random vector $(\boma{C}^t,\boma{R}^t)^t$.\\

\begin{theorem} 
\label{Teo:Assumption_Cov_matrices}
	If Assumption 1 holds then, for $j,l\in\{1,2,\ldots,p\}$ with $j\neq l$:
\begin{align}
	{\rm E}(A_j)&={\rm E}(C_j)\\
{\rm Var}(A_{j})&={\rm Var}(C_j)+\frac{1}{4}{\rm Var}(U_{j}){\rm E}(R_j^2)\label{Var(Aij)_Uj}	\\ 
{\rm Cov}(A_{j},A_{l})&={\rm Cov}(C_j,C_l)+ \frac{1}{4} {\rm Cov}(U_{j},U_{l}){\rm E}(R_jR_l). \label{Cov(Aij)_Uj}
\end{align}	
\end{theorem}

\begin{proof}
Let Assumption 1 hold. From the fact that ${\rm E}(U_{j})=0$,  $U_j$ is independent of $R_j$, and the linearity of the expectation operator, it follows that 
$${\rm E}(A_{j})={\rm E}(C_{j}) + \frac{1}{2} {\rm E}(U_{j}R_j)= {\rm E}(C_{j})+ \frac{1}{2} {\rm E}(U_{j})  {\rm E}(R_{j})={\rm E}(C_{j}).$$
%
% in order for micro-data and the associated macro-data having the same mean,  the weights $U_{j}$  should have zero mean. 

We will now obtain the symbolic variances and covariances. From \eqref{Aij} we can derive:
%In the first and less restricted case,  
\begin{align*}
{\rm Var}(A_{j})&={\rm Var}(C_j)+\frac{1}{4}{\rm Var}(U_{j}R_j)+{\rm Cov}(C_j,U_{j}R_j).
\end{align*}
Moreover, given that $U_{j}$ and $R_j$ are independent random variables and  ${\rm E}(U_{j})=0$, we obtain 
	$${\rm Var}(U_{j}R_j)={\rm E}(U_{j}^2){\rm E}(R_j^2)-{\rm E}(U_{j})^2{\rm E}(R_j)^2={\rm Var}(U_{j}){\rm E}(R_j^2).$$
Following a similar reasoning, it can be shown that 
	$${\rm Cov}(C_j,U_{j}R_j)={\rm E}(U_{j}){\rm E}(C_jR_j)-{\rm E}(C_j){\rm E}(U_{j}){\rm E}(R_j)=0$$
for all values of $j$. As a result, it follows as wanted that:
\begin{align*}
{\rm Var}(A_{j})&={\rm Var}(C_j)+\frac{1}{4}{\rm Var}(U_{j}){\rm E}(R_j^2).%\label{Var(Aij)_Uj}
\end{align*}

Again from \eqref{Aij} it follows that for $j,l\in\{1,2,\ldots,p\}$ with $j\neq l$;
\begin{align*}
{\rm Cov}(A_{j},A_{l})=&\, {\rm Cov}(C_j,C_l)+\frac{1}{2}{\rm Cov}(C_j,U_{l}R_l)+\frac{1}{2}{\rm Cov}(U_{j}R_j,C_l)+\\
&+\frac{1}{4}{\rm Cov}(U_{j}R_j,U_{l}R_l). 
\end{align*}
As $(U_{j},U_{l})$ is independent of $(C_j,R_j,C_l,R_l)$, it can be shown that 
\begin{eqnarray*}{\rm Cov}(U_{j}R_j,C_l)={\rm Cov}(C_j,U_{l}R_l)=0,\end{eqnarray*}
and 
\begin{eqnarray*}
	{\rm Cov}(U_{j}R_j,U_{l}R_l)&=&{\rm E}(U_{j}U_{l}){\rm E}(R_jR_l)-{\rm E}(U_{j}){\rm E}(U_{l}){\rm E}(R_j){\rm E}(R_l)\\
	&=&{\rm E}(U_{j}U_{l}){\rm E}(R_jR_l),
\end{eqnarray*}
 thus implying that:
\begin{align*}
	{\rm Cov}(A_{j},A_{l})&={\rm Cov}(C_j,C_l)+ \frac{1}{4} {\rm }{\rm Cov}(U_{j},U_{l}){\rm E}(R_jR_l). %\label{Cov(Aij)_Uj}
\end{align*}
\qedwhite %\smartqed 
%$	\smartqed $
\end{proof}

The following result, which provides the symbolic mean vector,   $\boma{\mu}={\rm E}(\boma{A})$, and the symbolic covariance matrix, 
$\boma{\Sigma} = {\rm Var}(\boma{A})$, follows trivially from Theorem~\ref{Teo:Assumption_Cov_matrices}. %%%%VERIFICAR APa

\begin{corollary}
 \label{Cor:symbmeancov}
 	If Assumption 1 holds then $\boma{\mu}={\rm E}(\boma{C})$ and
\begin{equation}
 \label{Cov(AU)}
	\boma{\Sigma}=  \boma{\Sigma} _{CC} 
	+ \frac{1}{4}  [\boma{\Sigma}_{UU} \bullet {\rm E}(\boma{R} \boma{R}^t)]
\end{equation}
where $\bullet$ denotes the Schur or entrywise product of matrices.
\end{corollary}

\subsection{Symbolic covariance and pseudo-correlation matrices under two scenarios}

We will now particularise the form of the symbolic covariance matrix for two scenarios more restrictive than Assumption 1 stated in the previous subsection: the first scenario in which the weights $U_j$ are uncorrelated random variables, and the second scenario in which the weights $U_j$ are almost surely equal (represented by $\overset{\rm a.s.}{=}$) to a random variable $U$. Adding these constraints to Assumption 1, namely that the weights $U_{j}$ are zero mean random variables with support on the interval $[-1,1]$,  independent from $(\boma{C}^t,\boma{R}^t)^t$, we have the following assumptions for scenarios 1 and 2:\\

\noindent {\bf Scenario 1}:  The weights $U_1,U_2,\ldots,U_p$ are zero mean uncorrelated random variables with support on $[-1,1]$ and are independent from $(\boma{C}^t,\boma{R}^t)^t$.\\

\noindent {\bf Scenario 2}:  $U_1 \overset{\rm a.s.}{=}U_2\overset{\rm a.s.}{=}\ldots \overset{\rm a.s.}{=} U_p\overset{\rm a.s.}{=} U$, with $U$ being a zero mean random variable with support on $[-1,1]$ and independent from $(\boma{C}^t,\boma{R}^t)^t$.

We note that Scenario 1 corresponds to the least possible linear association between the weights as they are uncorrelated random variables. In opposition, Scenario 2 corresponds to the highest possible association between the weights as they are equal (almost surely).

The result  \eqref{Cov(AU)}, along with Scenarios 1 and 2,  lead to two different families of symbolic covariance matrices, summarised in Corollary~\ref{Cor:Families_of_Cov_matrices}. In what follows, if $\boma{Z}$ is a matrix then ${\rm Diag} (\boma{Z})$ represents a diagonal matrix whose main diagonal have the values $[\boma{Z}]_{ii}$. If $\boma{z}$ is a vector then ${\rm Diag} (\boma{z})$ represents a diagonal matrix whose main diagonal have the values of vector  $\boma{z}$.

\begin{corollary} \label{Cor:Families_of_Cov_matrices}
	Under Scenario 1 the symbolic covariance matrix, denoted by $\boma{\Sigma}^{(1)}$, has the form
				\begin{equation}
				\label{Eq:SymCov_Diag_DifVar(U_j)}
				 \boldsymbol \Upsigma^{(1)} = \boldsymbol \Upsigma_{CC} + \displaystyle \frac{\boma{\Sigma}_{UU}}{4}\;{\rm Diag} \left({\rm E}(\boma{RR}^t)\right),
				\end{equation}
				where $\boma{\Delta}={\rm Diag} \left({\rm Var}(U_1),\ldots,{\rm Var}(U_p)\right)$. Moreover, if 
				all $U_j$ share the same variance (i.e. ${\rm Var}(U_j)={\rm Var}(U_1)$ for $j\in\{2,3,\ldots,p\}$), then $\boma{\Sigma}^{(1)}$ simplifies to:
				\begin{equation}
					\label{Eq:SymCov_Diag}
				\boldsymbol \Upsigma^{(1)} = \boldsymbol \Upsigma_{CC} + \displaystyle \frac{{\rm Var}(U_1)}{4}\;{\rm Diag} \left({\rm E}(\boma{RR}^t)\right).
				\end{equation}	
				
	Under Scenario 2 the symbolic covariance matrix, denoted by $\boma{\Sigma}^{(2)}$,  has the form							
				\begin{equation}\label{Eq:SymCov_S2}
				\boldsymbol \Upsigma^{(2)} = \boldsymbol \Upsigma_{CC} + \displaystyle \frac{{\rm Var}(U)}{4}\;{\rm E}(\boma{RR}^t).
				\end{equation}
\end{corollary}

\begin{proof}
Let us first assume that the conditions of Scenario 1 hold. This implies that the conditions of Assumption 1 hold as well and, in particular, that equation (\ref{Cov(AU)})
holds. As, in addition, $\boma{\Sigma}_{UU}={\rm Diag} \left({\rm Var}(U_1),\ldots,{\rm Var}(U_p)\right)$, since the weights $U_1$, $U_2$, \ldots, $U_p$ are uncorrelated random variables, we conclude that  (\ref{Eq:SymCov_Diag_DifVar(U_j)})  holds. Moreover, as (\ref{Eq:SymCov_Diag}) follows trivially from  (\ref{Eq:SymCov_Diag_DifVar(U_j)}), we conclude that  (\ref{Eq:SymCov_Diag}) holds as well.

Let us now assume that the conditions of Scenario 2 hold. This implies that the conditions of Assumption 1 hold as well and, in particular, that equation (\ref{Cov(AU)})
holds. As, in addition,   $ {\rm  Cov}(U_{j},U_{l})= {\rm E}(U_{j}U_{l})= {\rm E}(U^2)= {\rm Var}(U)$, for $j,l\in\{1,2,\ldots,p\}$ since 
$U_1 \overset{\rm a.s.}{=}U_2\overset{\rm a.s.}{=}\ldots \overset{\rm a.s.}{=} U_p\overset{\rm a.s.}{=} U$,
we conclude that  (\ref{Eq:SymCov_S2})  holds. 
\qedwhite 
\end{proof}

Different particular forms of the symbolic covariance matrices $\boldsymbol \Upsigma^{(1)}$ and $\boldsymbol \Upsigma^{(2)}$ will be derived in the next subsection, along with conditions on  the weights $U_j$ that may lead to such symbolic covariance matrices. However, we will first address  pseudo-correlation matrices associated to the matrices  $\boldsymbol \Upsigma^{(1)}$ and $\boldsymbol \Upsigma^{(2)}$, starting by the definition of pseudo-correlation matrix associated with a square matrix with positive diagonal entries.

 \begin{definition}\label{def:PseudoCorr}
		If $\boma{\Lambda}=[\lambda_{jl}]_{j,l\in\{1,2,\ldots,p\}}$ is a square matrix with positive diagonal entries, then the associated pseudo-correlation matrix, 
	denoted by $\rho( \boma{\Lambda})=[[\rho(\boma{\Lambda})]_{jl}]_{j,l\in\{1,2,\ldots,p\}}$, is given by
	\begin{equation} \label{def:PseudoCorrjl}
		[\rho(\boma{\Lambda})]_{jl} = \frac{\lambda_{jl}}{\sqrt{\lambda_{jj} \lambda_{ll}}},\quad  j,l\in\{1,2,\ldots,p\}.
	\end{equation}	
\end{definition}

Note that, in particular, $\rho( \boma{\Lambda})$ is the correlation matrix associated to a given $p$-dimensional random vector case $\boma{\Lambda}$ is the covariance matrix of the same random vector. 

In the next theorem we show that the absolute values of the entries of the pseudo-correlation matrices associated to the matrices $\boldsymbol \Upsigma^{(1)}$ and $\boldsymbol \Upsigma^{(2)}$  (see equations (\ref{Eq:SymCov_Diag_DifVar(U_j)}) and (\ref{Eq:SymCov_S2}))  are smaller or equal to one. In other words, 
symbolic pseudo-correlations associated to the matrices  $\boldsymbol \Upsigma^{(1)}$ and $\boldsymbol \Upsigma^{(2)}$ are quantities between -1 and 1, as in the conventional case [provided $\boldsymbol \Upsigma^{(1)}$ and $\boldsymbol \Upsigma^{(2)}$ have finite entries, which is essentially the case whenever $C_1,C_2,\ldots,C_p,R_1,R_2,\ldots,R_p$ are random variables with finite moments of order two].

\begin{theorem} \label{Desig_Cauchy_Schwarz_SymCor}
	If $\boma{X}$ is such that  $\boldsymbol \Upsigma^{(1)}$ in the form given by equation (\ref{Eq:SymCov_Diag})  has only positive finite entries, then 
\begin{equation}
\label{ineq:Cauchy_Schwarz_SymCorr1}
	\left|[\rho(\boldsymbol \Upsigma^{(1)})]_{jl}\right| \leq 1
	\quad\quad \text{for} \quad\quad j,l\in\{1,2,\ldots,p\}.
\end{equation}

Similarly, if $\boma{X}$ is such that  $\boldsymbol \Upsigma^{(2)}$ given by equation (\ref{Eq:SymCov_S2}) has only positive finite entries, then 
\begin{equation}
\label{ineq:Cauchy_Schwarz_SymCorr2}
	\left|[\rho(\boldsymbol \Upsigma^{(2)})]_{jl}\right| \leq 1
	\quad\quad \text{for} \quad\quad j,l\in\{1,2,\ldots,p\}.
\end{equation}
\end{theorem} 

\begin{proof}
	We first suppose that  $\boldsymbol \Upsigma^{(1)}$ given by equation (\ref{Eq:SymCov_Diag}) has only finite positive entries, in which case 
$[\rho(\boldsymbol \Upsigma^{(1)})]_{jj}=1$, for all $j\in\{1,2,\ldots,p\}$, follows trivially. Secondly, we consider arbitrary indices $j,l\in\{1,2,\ldots,p\}$ such that $j\neq l$ 
and note that, using the Cauchy-Schwarz inequality applied to the random pair $(C_j,C_l)$, we obtain, as wanted,
	\begin{eqnarray*}
		\left|[\rho(\boldsymbol \Upsigma^{(1)})]_{jl}\right| &= &
		\frac{\left|{\rm Cov}(C_j,C_l)\right|}{\sqrt{[{\rm Var}(C_j)+{\rm Var}(U_1)\text{E}(R_j^2)]  [{\rm Var}(C_l)+{\rm Var}(U_1)\text{E}(R_l^2)]}}\\
		&\leq &\frac{\sqrt{{\rm Var}(C_1){\rm Var}(C_2)}}
			{\sqrt{[{\rm Var}(C_j)+{\rm Var}(U_1)\text{E}(R_j^2)]  [{\rm Var}(C_l)+{\rm Var}(U_1)\text{E}(R_l^2)]}}\\
		&\leq &1.
	\end{eqnarray*}	

We now assume that  $\boldsymbol \Upsigma^{(2)}$ given by equation (\ref{Eq:SymCov_S2}) has only positive finite entries, in which case 
$[\rho(\boldsymbol \Upsigma^{(2)})]_{jj}=1$, for all $j\in\{1,2,\ldots,p\}$, follows trivially. We next consider arbitrary indices $j,l\in\{1,2,\ldots,p\}$ such that $j\neq l$ 
and note that, using the Cauchy-Schwarz inequality applied to the random pair $(C_j,C_l)$, Holder's inequality applied to te random pair $(R_j,R_l)$,  and making $\delta = {\rm Var}(U)/4$, we obtain:
	\begin{eqnarray*}
			\left|[\rho(\boldsymbol \Upsigma^{(2)})]_{jl}\right| &\leq &
			\frac{|{\rm Cov}(C_j,C_l)|+\delta\, {\rm E}(R_jR_l)}
			{\sqrt{[{\rm Var}(C_j)|+\delta {\rm E}(R_j^2]  [{\rm Var}(C_j)|+\delta\, {\rm E}(R_j^2] }} \\
			&\leq & 	\frac{\sqrt{{\rm Var}(C_j)  {\rm Var}(C_l)} + \delta \sqrt{ {\rm E}(R_j^2)  {\rm E}(R_l^2 )}}
			{\sqrt{[{\rm Var}(C_j)|+\delta\, {\rm E}(R_j^2]  [{\rm Var}(C_j)|+\delta\, {\rm E}(R_j^2] }}.%\label{SCor_Desig1}
	\end{eqnarray*}
Thus, $\left[[\rho(\boldsymbol \Upsigma^{(2)})]_{jl}\right]^2$ is smaller or equal to
	\begin{align*}
	\frac{{\rm Var}(C_j){\rm Var}(C_l)+\delta^2 {\rm E}(R_j^2){\rm E}(R_l^2)+2\delta {\sqrt{{\rm Var}( C_j){\rm E}(R_l^2){\rm Var}( C_l){\rm E}(R_j^2)}}}
	%Denominator
	{{\rm Var}(C_j){\rm Var}(C_l)+\delta^2 {\rm E}(R_j^2){\rm E}(R_l^2)+\delta \left({\rm Var}(C_j){\rm E}(R_l^2)+{\rm Var}(C_l){\rm E}(R_j^2)\right)}.
	\end{align*}
As it can be easily proved, for all non-negative values of $x$ and $y$, $2\sqrt{xy} \leq (x+y)$. Thus, making $x={\rm Var}(C_j){\rm E}(R_l^2)$ and 
$y={\rm Var}(C_l){\rm E}(R_j^2)$, we conclude that, as wanted,  $\left[[\rho(\boldsymbol \Upsigma^{(2)})]_{jl}\right]^2 \leq 1$.  \qedwhite
\end{proof}

\subsection{Eight particular forms of symbolic covariance matrices}
\label{8SymbCovMatrices}

\begin{figure}
\centering
\captionsetup{justification=raggedright}
  \subfloat[$\boma{\Sigma}=\boma{\Sigma}_1$.]{\label{fig:sigma1}\includegraphics[width=.3\textwidth]{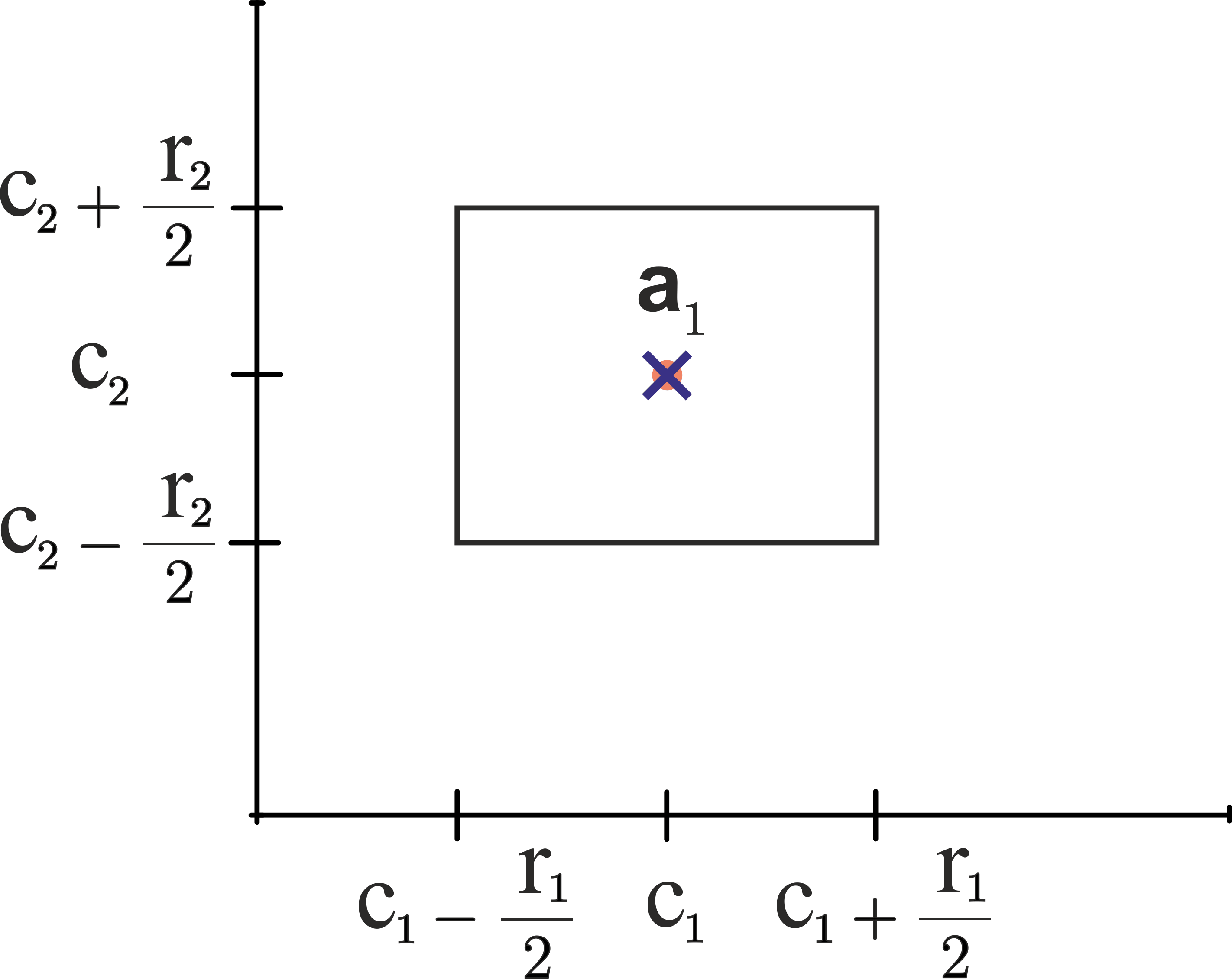}}                
  \hspace*{5pt}
  \subfloat[$\boma{\Sigma}=\boma{\Sigma}_2$.]{\label{fig:sigma2}\includegraphics[width=.3\textwidth]{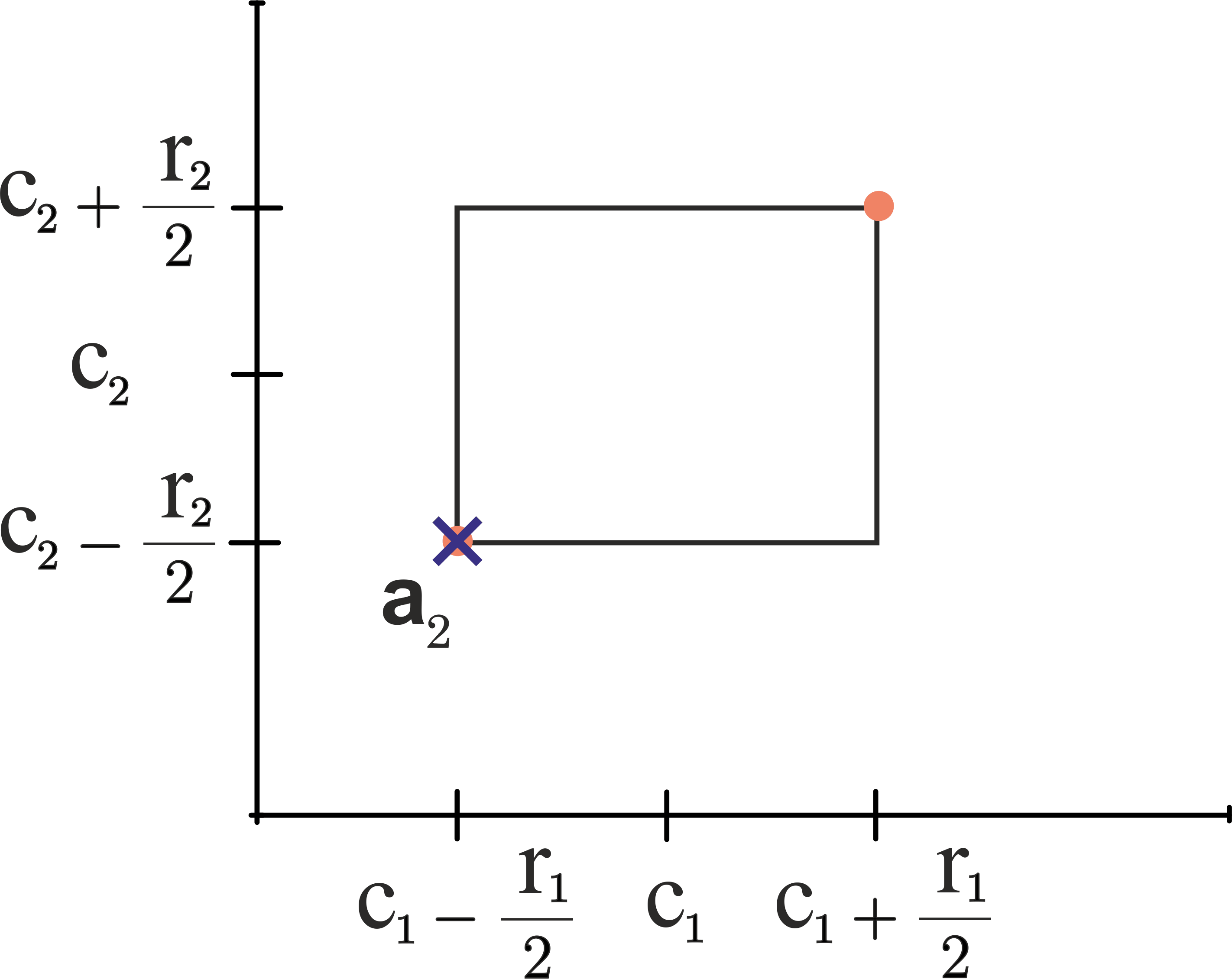}}
    \hspace*{5pt}
  \subfloat[$\boma{\Sigma}=\boma{\Sigma}_3$.]{\label{fig:sigma3}\includegraphics[width=.3\textwidth]{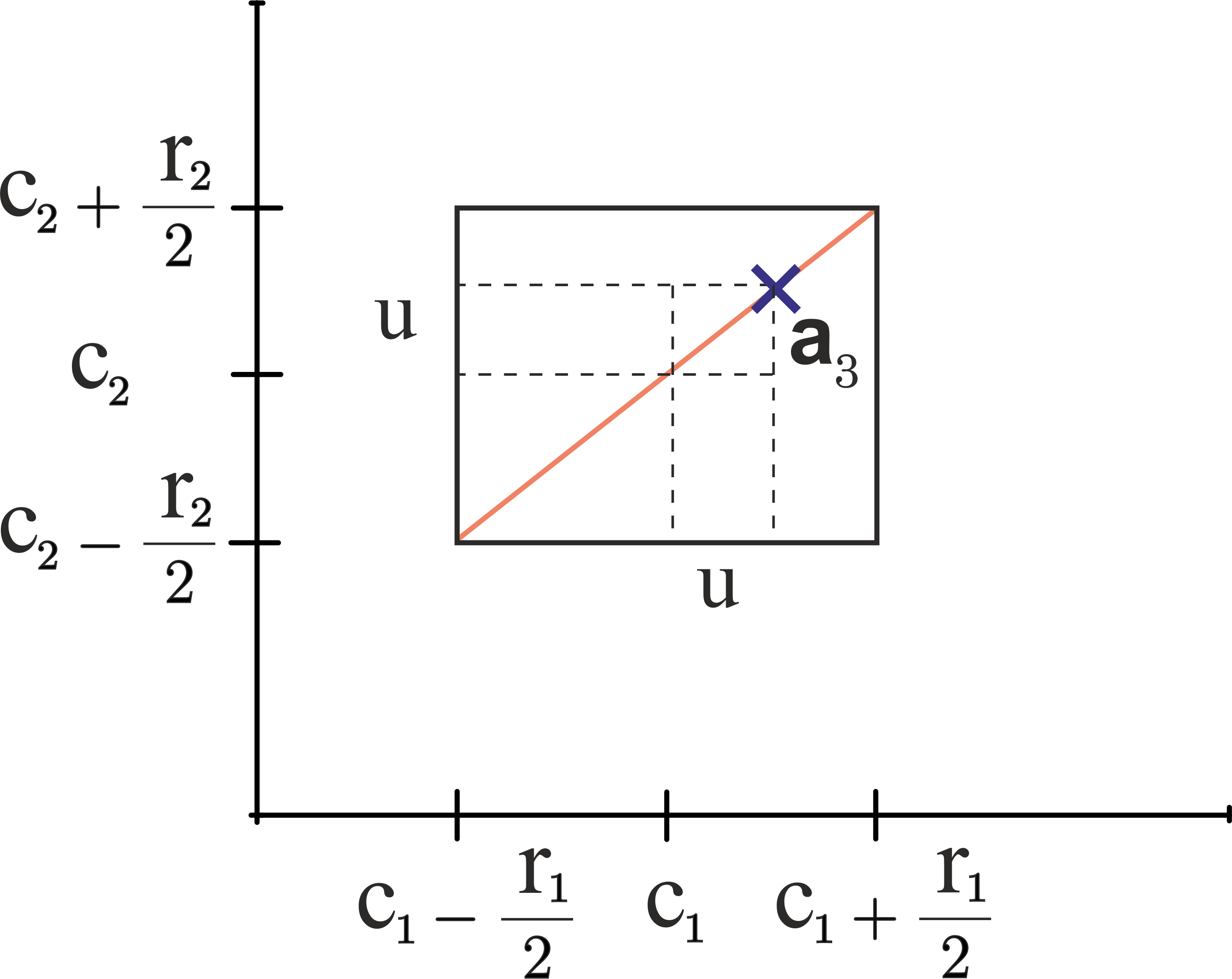}}                \\
  \subfloat[$\boma{\Sigma}=\boma{\Sigma}_4$.]{\label{fig:sigma4}\includegraphics[width=.3\textwidth]{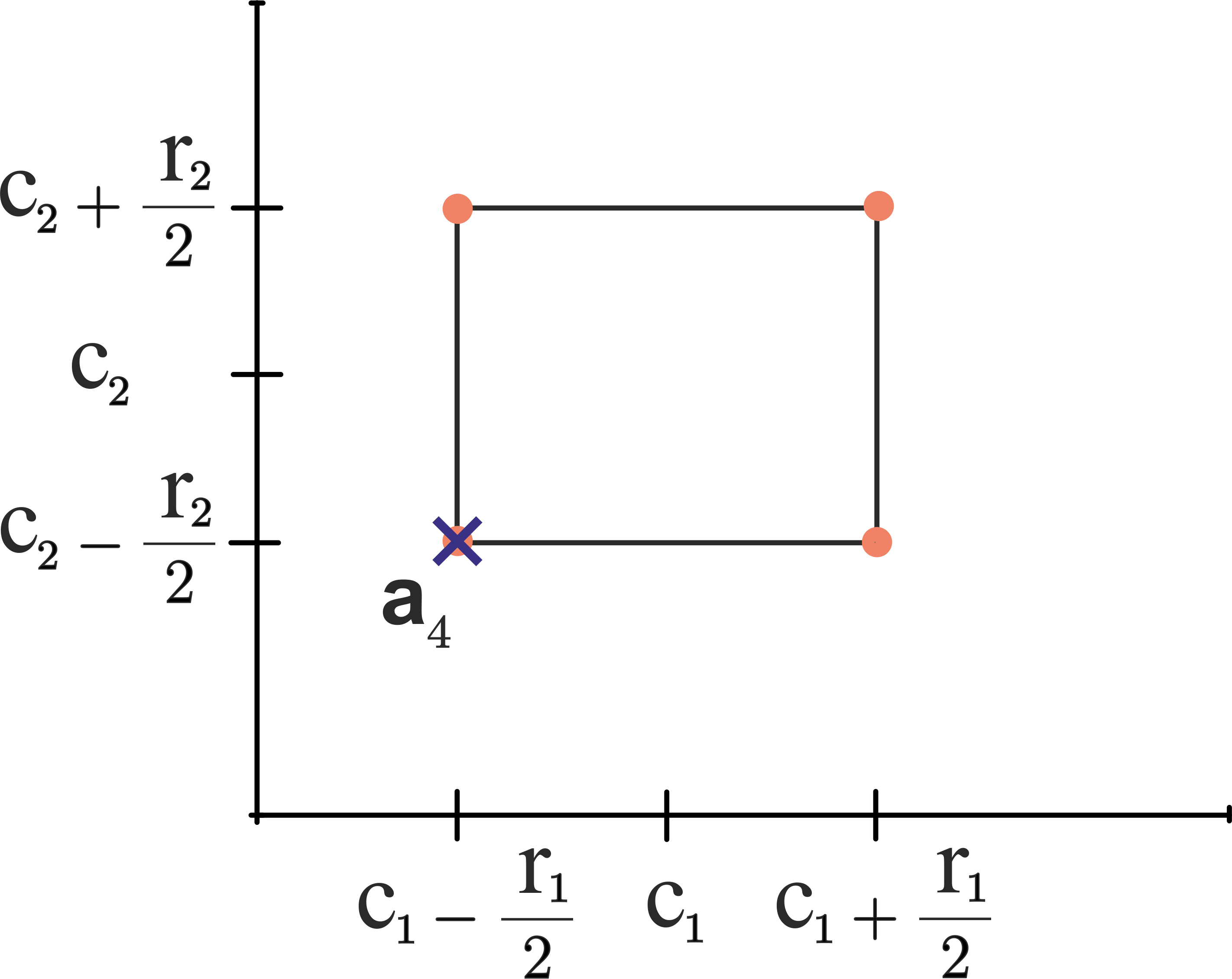}}
    \hspace*{5pt}
  \subfloat[$\boma{\Sigma}=\boma{\Sigma}_k,\, k=5,\ldots,8$.]{\label{fig:sigma5}\includegraphics[width=.3\textwidth]{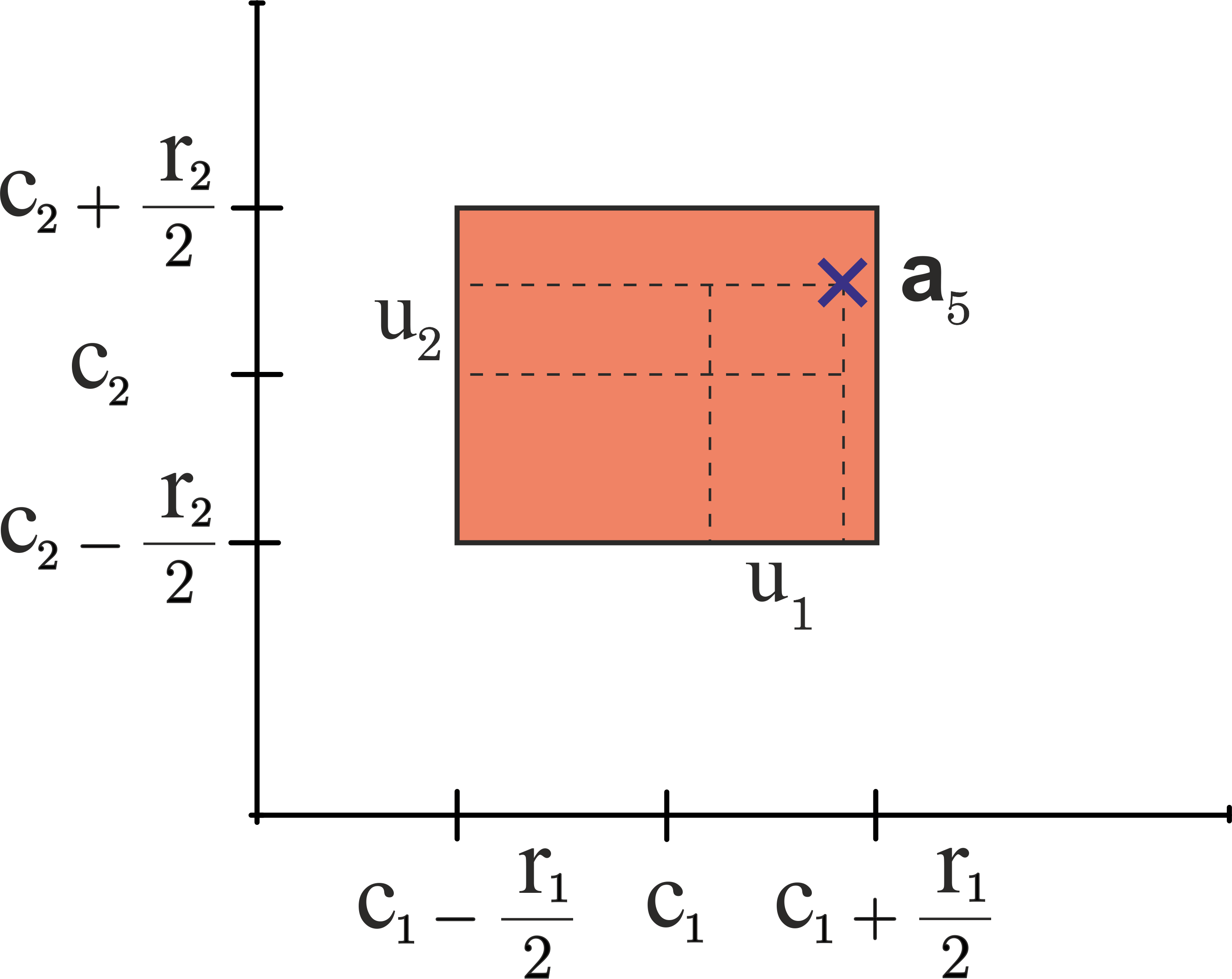}}          
\caption{Representation of possible values of $\boma{A}=(A_{1},A_{2})^t$, such that $\boma{\Sigma}=\boma{\Sigma}_k$, $k=1,\ldots,8$.}
  \label{fig:ReasoningSigmaFormulation}
\end{figure}

In Table~\ref{tab:reasoningVarCov} we list 8 particular forms for the symbolic covariance matrix, $\boma{\Sigma}$, along with the scenario and additional constraints on the weights $U_1, U_2,\ldots, U_p$ or $U$ that may give rise to them.  Moreover, in Figure \ref{fig:ReasoningSigmaFormulation}, we provide an illustration of the associated micro-data for the 2-dimensional case; i.e., we display the possible values of $\boma{A}=(A_{1},A_{2})^t$, for $k=1,2,\ldots,8$.%RV 3/11/2017

	\begin{table}[H]
		\caption{List of  symbolic covariance matrices $\boma{\Sigma}_k$, $k=1,2,\ldots,8$, and scenario along with additional constraints on the weights $U_j$ or $U$ under which  one has  $\boma{\Sigma}={\rm Cov}(\boma{A})=\boma{\Sigma}_k$.}
		\label{tab:reasoningVarCov}
		\centering
		\scriptsize
		\setlength{\tabcolsep}{2.5pt}
		\begin{tabular}{clcl}
			\toprule
			\multicolumn{1}{c}{\begin{tabular}{c}
		Symbolic\\ Covariance\\ Matrix
			\end{tabular}} &  \multicolumn{1}{c}{Definition} &  Scenario & \multicolumn{1}{c}{$(U_1, U_2,\ldots, U_p)$  or $U$} \\
			\midrule
			$ \boldsymbol \Sigma_1$  & $\boldsymbol \Sigma_{CC}$ &  1 or 2 & $P(U_{j}=0)=1,\,\forall j$ or $P(U=0)=1$\\[8pt]
			$ \boldsymbol \Sigma_2$ & $\boldsymbol \Sigma_{CC} +  \displaystyle \frac{1}{4}\; {\rm E}(\boma{RR}^t)$ &  2 & $U \sim {\rm Unif}\{-1,1\}$  \\[8pt]
			 $ \boldsymbol \Sigma_3$& $\boldsymbol \Upsigma_{CC} + \displaystyle \frac{1}{12}\;  {\rm E}(\boma{RR}^t)$ & 2 & ${\rm Var}(U)=\frac{1}{3}$ \\[8pt]
			 $ \boldsymbol \Sigma_4$& $\boldsymbol \Upsigma_{CC} + \displaystyle \frac{1}{4}\;{\rm Diag} \left({\rm E}(\boma{RR}^t)\right)$ &  1 & $U_j \sim {\rm Unif}\{-1,1\},\, \forall j$ \\[8pt]
			 $ \boldsymbol \Sigma_5$&  $\boldsymbol \Upsigma_{CC} + \displaystyle \frac{1}{12}\;{\rm Diag} \left({\rm E}(\boma{RR}^t)\right)$ & 1 & ${\rm Var}(U_{j})=\frac{1}{3},\, \forall j$ \\[8pt]
			$ \boldsymbol \Sigma_6$& $\boldsymbol \Upsigma_{CC} + \displaystyle \frac{1}{8}\;{\rm Diag} \left({\rm E}(\boma{RR}^t)\right)$ &  1 &$U_{j} \sim {\rm  InvTriangular}[-1,0,1]$\\[8pt] %& $|u|$\\[8pt]
			$ \boldsymbol \Sigma_7$& $\boldsymbol \Upsigma_{CC} + \displaystyle \frac{1}{24}\;{\rm Diag} \left({\rm E}(\boma{RR}^t)\right)$ &  1 &$U_{j} \sim {\rm  Triangular[}-1,0,1]$\\[8pt] %& $1-|u|$\\[8pt]
			$ \boldsymbol \Sigma_8$&  $\simeq \boldsymbol \Upsigma_{CC} + \displaystyle \frac{1}{36}\;{\rm Diag} \left({\rm E}(\boma{RR}^t)\right)$\footnote{The exact weight of the ranges contribution is: $\delta_8=\frac{1}{36}-\frac{ \phi(3)}{6(2\Phi(3)-1)}$.} & 1 &$U_{j}=Z|Z\in [-1,1]$, $Z \sim {\cal N}(0,\frac{1}{9})$ \\[8pt] %& $\frac{\dps 3\phi( 3u)}{\dps  (2\Phi(3)-1)}$\\[8pt]

			\bottomrule
		\end{tabular}
	\end{table}

\captionsetup[subfigure]{subrefformat=simple,labelformat=simple,listofformat=subsimple}
\renewcommand\thesubfigure{(\arabic{subfigure})}

The first form ($k=1$) is obtained under Scenario 1 with the weights $U_{j}$ being zero with probability one, or equivalently  under Scenario 2 with the weight $U$ being zero with probability one, corresponding to the case where the ranges are not taken into account. In such case $\boma{A} \overset{\rm a.s.}{=}\boma{C}$, leading to $ \boldsymbol \Sigma = \boldsymbol \Sigma_1={\rm Cov}(\boma{C})$. As illustrated in Figure \subref*{fig:sigma1} for the 2-dimensional case, the micro-data associated with this model is always at the center of the interval-valued object.%RV 3/11/2017

The second form ($k=2$) is obtained under Scenario 2 when the weight $U$ has unit variance; this along with $U$ having zero mean and support on $[-1,1]$ implies that $U$ follows a discrete Uniform distribution on $\{-1,1\}$. Then, ${\rm Var}(U)=1$ leads to $\boma{\Sigma}= {\rm Cov}(\boma{A})=\boma{\Sigma}_2$, in view of \eqref{Eq:SymCov_S2}. The micro-data associated with this model is at one of the two vertices  $(C_1-R_1/2,C_2-R_2/2,\ldots,C_p-R_p/2)^t$ or $(C_1+R_1/2,C_2+R_2/2,\ldots,C_p+R_p/2)^t$, chosen with equal probability. The 2-dimensional case is illustrated in Figure~\subref*{fig:sigma2}.%RV 3/11/2017

The third form ($k=3$) is obtained under Scenario 2 when the weight $U$ has variance $1/3$, as it is the case when $U$ follows a continuous Uniform distribution on $[-1,1]$. Then, ${\rm Var}(U)=1/3$  leads to $\boma{\Sigma}={\rm Cov}(\boma{A})=\boma{\Sigma}_3$, in view of  \eqref{Eq:SymCov_S2}. Note that when $U$ follows a continuous Uniform distribution on $[-1,1]$, the micro-data associated with the model is at the line segment with endpoints $(C_1-R_1/2,C_2-R_2/2,\ldots,C_p-R_p/2)^t$ and $(C_1+R_1/2,C_2+R_2/2,\ldots,C_p+R_p/2)^t$, which passes through the vector of centers 
$(C_1,C_2,\ldots,C_p)^t$. Figure~\subref*{fig:sigma3} illustrates the 2-dimensional case.%RV 3/11/2017

With a similar reasoning, we may conclude that the fourth form ($k=4$) is obtained under Scenario 1 when the weights $U_{1},U_2,\ldots,U_{p}$ are independent random variables with discrete Uniform distribution on $\{-1,1\}$. Since ${\rm Var}(U_j)=1$, for all $j$, in view of  \eqref{Eq:SymCov_Diag}, we conclude that $\boma{\Sigma}={\rm Cov}(\boma{A})=\boma{\Sigma}_4$. Micro-data associated with this model are at one of the vertices of the hyper-rectangle corresponding to the Cartesian product $(C_1-R_1/2,C_1+R_1/2)\times (C_2-R_2/2,C_2+R_2/2)\times \ldots \times (C_p-R_p/2,C_p+R_p/2)$, with  the same probability, $1/2^p$, as illustrated in Figure \subref*{fig:sigma4}, for the 2-dimensional case.%RV 3/11/2017

The fifth form ($k=5$) is obtained under Scenario 1 when $U_{1},U_2,\ldots,U_{p}$ are uncorrelated random variables with common variance $1/3$; note that is the case when the weights $U_j$ are  independent random variables with continuous uniform distribution on $[-1,1]$, independent from $(\boma{C},\boma{R})$. Since ${\rm Var}(U_j)=1/3$, for all $j$, in view of \eqref{Eq:SymCov_Diag}, it follows that $\boma{\Sigma}={\rm Cov}(\boma{A})=\boma{\Sigma}_5$. In case  the weights $U_j$ are  independent random variables with continuous uniform distribution on $[-1,1]$, independent from 
$(\boma{C}^t,\boma{R}^t)^t$, the micro-data follow a  $p$-dimensional continuous uniform distribution on the hyper-rectangle corresponding to the Cartesian product $(C_1-R_1/2,C_1+R_1/2)\times (C_2-R_2/2,C_2+R_2/2)\times \ldots \times (C_p-R_p/2,C_p+R_p/2)$. This case is illustrated in Figure \subref*{fig:sigma5}, for the 2-dimensional case. 

The sixth form ($k=6$) models micro-data that can take any value inside the immediately above mentioned hyper-rectangle but showing a tendency to be located near its vertices. The points are spread according to a $p$-dimensional inverse triangular distribution. 
Finally, the remaining two forms ($k=7$ and $k=8$) model micro-data mainly located near the center of the same hyper-rectangle, with more concentration near the center of the hyper-rectangle in the case $k=8$ ($p$-dimensional multivariate truncated normal distribution) than in the case $k=7$ ($p$-dimensional multivariate triangular distribution). %RV 24/11/2017
The covariance matrices $\boma{\Sigma}_6$ and $\boma{\Sigma}_7$ can be obtained from \eqref{Eq:SymCov_Diag} noting that: 
${\rm Var}(U_j)=1/2$ for $k=6$, and ${\rm Var}(U_j)=1/6$ for $k=7$. In case $k=8$, the form of $\boma{\Sigma}_8$ follows from the fact that if $Z$ has normal distribution with zero mean and ${\rm Var}(Z)=1/9$, $Z \sim {\cal N}(0,\frac{1}{9})$, then ${\rm Var}(U_j)={\rm Var}(Z|Z\in [-1,1])=\frac{1}{9}-\frac{2 \phi(3)}{3(2\Phi(3)-1)}$, where $\frac{2 \phi(3)}{3(2\Phi(3)-1)}\simeq 2.96\times 10^{-3}$, and $\phi(\cdot)$ 
and $\Phi(\cdot)$ are, respectively, the probability density function and cumulative distribution function of the standard normal distribution. 
%The covariances and micro-data models corresponding to cases $k=6,7,8$ are summarized in Table \ref{tab:reasoningVarCov}. %\ref{tab:reasoningVarCov_Cazes} 
% Note that these covariance matrices inherit all the properties of cases $k=4,5$ (see Proposition \ref{Prop_SymCov_I}), since the $U_{j}$ are independent, identically distributed symmetric random variables, with support $[-1,1]$, independent from 
%$(\boma{C},\boma{R})$.%RV 3/11/2017
%
%Among the eight covariance forms displayed in Table~\ref{tab:reasoningVarCov}, $\boma{\Sigma}_5$, $\boma{\Sigma}_6$, $\boma{\Sigma}_7$ and $\boma{\Sigma}_8$
%seem to be the ones with wider applicability, given the way the micro-data is allowed to be distributed inside the hyper-rectangle corresponding to the Cartesian product $(C_1-R_1/2,C_1+R_1/2)\times (C_2-R_2/2,C_2+R_2/2)\times \ldots \times (C_p-R_p/2,C_p+R_p/2)$.

Cazes et al. \cite{Cazes1997} have also addressed the modelling of micro-data, but under the restriction of fixed macro-data (i.e. macro-data with deterministic interval limits). In their work, they have considered micro-data models that fit under Scenario 1, as it is the case of the forms $k=1,4,5,6,7,8$. However, they have not related the structure of micro-data with the definitions of symbolic covariances matrices for interval-valued random variables. Specifically, they considered $U_{1},U_2,\ldots,U_{p}$ to be independent random variables with (i) inverse triangular distribution with parameters $\{-1,0,1\}$ (we refer to this case as $k=6$), (ii) triangular distribution with parameters $\{-1,0,1\}$ ($k=7$), and (iii)  Normal distribution truncated to $[-1,1]$, with zero mean and standard deviation approximately equal to 1/3 ($k=8$). These constraints on the micro-data  are worth considering since they bring alternative modelling possibilities to the forms $k=1,4,5$.%RV 3/11/2017

\subsection{How to choose a particular form of symbolic covariance matrix?}

The formulation \eqref{Aij} allows us to write:
\begin{equation}\label{Ukj}
U_{j}=\frac{2(A_{j}-C_{j})}{R_{j}}
\end{equation} 
for $j=1,2,\ldots,p$. 
Thus, if we have a dataset for which macro-data and micro-data are available, we can obtain 
a realization  $u_{j}$ of $U$ using the relation $u_j=2(a_j-c_j)/r_j$, with $a_j$, $c_j$, and $r_j$ being the realizations of $A_j$, $C_j$, and $R_j$, respectively. 
Based on $u_{j}$'s, evidences about the distributional form of the random weights $U_{j}$ may be explored, leading to the choice of the symbolic covariance (correlation) matrix (see Table \ref{tab:reasoningVarCov}) that better fits the data under study.%RV 3/11/2017

\begin{figure}
\centering
\includegraphics[width=\textwidth]{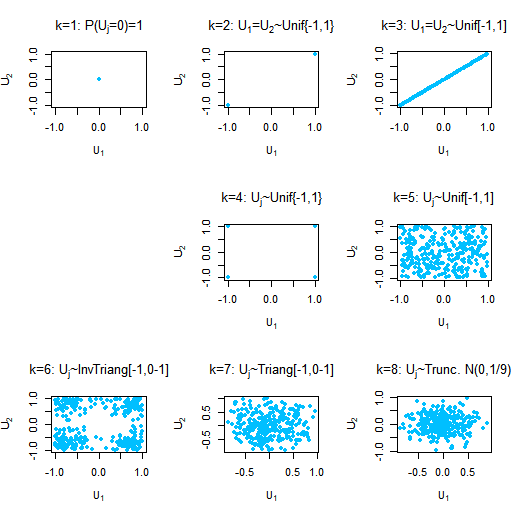}
\caption{Simulated values of $(U_{1},U_{2})^t$ according to the symbolic models listed in Table \ref{tab:reasoningVarCov} for  $k=1,\ldots,8$.}
  \label{fig:Simul_Aij}
\end{figure}

The scatter plots in Figure \ref{fig:Simul_Aij}, represent 300 simulated values of $(U_{1},U_{2})^t$, according to the eight models described in Table \ref{tab:reasoningVarCov}, illustrating typical patterns of points following each of the distributional models considered in the previous subsection, $k=1,2,\ldots,8$. 

Note that we may use the procedure described in this section for other models not considered above under Scenario 1 and Scenario 2. The practitioner may obtain new symbolic covariance matrices (based on new distributional forms for the random weights $U_{j}$) better suited to particular datasets.

In subsections \ref{Sec:Iris} and \ref{Sec:CreditCard}, we exemplify the choice of the distributional model for the $U_j$'s based on the available micro-data and macro-data for two specific datasets.

%% file: Rev_State_of_the_Art_v2.tex
\subsection{State of the art}\label{sec:SymCov}

In the literature, there are several alternative definitions of sample symbolic covariance matrix. In all definitions, the sample symbolic mean vector is $\bar{\boma{x}}=(\bar x_1,\ldots,\bar x_p)^t$ where
%\begin{align}
$\overline{x}_j  = \sum_{i=1}^{n} c_{ij},$ $j=1,\ldots,p$.    
%\end{align}
Choosing the sample mean of the centers as the sample symbolic mean makes sense, in particular, under the assumption that the micro-data associated with an interval follows a symmetric distribution on that interval.%RV 2/11/2020

In Oliveira et al. \cite{Oliveira.et.al:2017} the sample counterparts of definitions 1 to 5 of covariance matrices were written in terms of the sample covariance matrix of the centers and the sample second order moment of the ranges. The population versions of the existing proposals are listed in Table \ref{tab:reasoningVarCov}.%RV 2/11/2020 

The sample equivalent of $\boma{\Sigma}_1$ was proposed by Billard and Diday \cite{BillardDiday2003}, and is the most straightforward approach, since in its essence follows the definition of the conventional sample variance and covariance of the interval centers. It has the disadvantage of ignoring the  contribution of the ranges to the sample symbolic variance and covariance. To overcome this limitation, de Carvalho et al. \cite{Carvalho.et.al:2006} proposed a sample variance definition based on the squared distances between the interval limits, and generalized this idea to define the associated sample covariance. Such definitions led to the sample equivalent of $\boma{\Sigma}_2$. A third alternative, proposed by Bertrand and Goupil \cite{Bertrand2000}, is obtained from the empirical density function of an interval-valued variable, assuming that the micro-data follows a Uniform distribution. The corresponding covariance was introduced by  Billard \cite{Billard2008}, and is based on the explicit decomposition of the covariance into within sum of products and between sum of products. Jointly they define the sample equivalent of $\boma{\Sigma}_3$.%RV 2/11/2020 

Other sample alternatives for covariance matrices were proposed in the context of symbolic principal component analysis (SPCA) for interval data. In one of the SPCA approaches, the authors define a symbolic sample covariance and use their eigenvectors as the weights of the sample symbolic principal variables. The sample version of $\boma{\Sigma}_4$ is part of the SPCA method called the {\em vertices method} and introduced in \cite{Cazes1997}. In the same paper, the authors considered other distribution for the micro-data that led to the  sample versions of  definitions 6 to 8 of Table  \ref{tab:reasoningVarCov}. The sample counterpart of $\boma{\Sigma}_5$ was also proposed as a part of a symbolic principal component estimation method called Complete-Information-based Principal Component \cite{Wang2012}.%RV 2/11/2020 

%% file: Rev_Null_SCov_v2.tex
\subsection{Null Symbolic Correlation}\label{sec:NullSymCov}

One of the reasons that makes the correlation coefficient so popular and used is its ability to quantify the strength of linear associations between pairs of random variables. Modeling real problems using linear combinations has the appeal of leading to understandable interpretations and the associated mathematical problems are, in general, easily solved. Linear combinations of input variables have been of central importance in  classical multivariate methodologies, like principal components analysis, factor analysis, canonical correlation analysis, discriminant analysis, linear regression models, among others \cite{Johnson2007}. In addition, its importance seems to have been reinvented with deep learning approaches \cite{Lecun.et.al:2015}.%RV 1/11/2020

Given that the correlation is a measure of linear association, a null value does not ensure the independence between random variables, but gives indication that there is no linear (and only linear) association between them. In the symbolic framework, the meaning of uncorrelated interval-valued variables is of importance and is analyzed in this subsection.%RV 1/11/2020

As discussed previously, in the literature there are two families of definitions of symbolic correlations (covariances). In the first family (Scenario 1), only the association between the centers of the interval-valued random variables are taken into consideration, namely ${\rm Cov}_k(X_j,X_l)={\rm Cov}(C_j,C_l)$, for $k=1,4,5,\ldots,8$ (see Table~\ref{tab:reasoningVarCov}); to make the discussion easier, we only consider cases $k=1,4,5$ of this family. The second family (Scenario 2) takes the ranges' contribution into consideration by adding the second order cross-moment of the ranges, weighted according to each definition: ${\rm Cov}_k(X_j,X_l)={\rm Cov}(C_j,C_l)+\delta_k{\rm E}(R_jR_l)$, for $k=2,3$, with $\delta_2=1/4$ and $\delta_3=1/12$. The first observation is that any association between centers and ranges that may exist is not detected by the symbolic covariances, since the definitions do not include any cross term between centers and ranges. Another potential pitfall is that we can devise cases where existing associations among ranges or among centers are not detected by the symbolic covariance definitions. Table \ref{tab:NullPitfalls} summarizes some of these problems, which will be discussed next.%RV 1/11/2020

\begin{minipage}[t]{0.9\textwidth}%{\linewidth \labelsep 1pt}
	\begin{table}[H]
	\caption{Cases where the  symbolic covariance definitions lead to unexpected results.}
	\small
	\begin{tabular}{lclcc}
		\hline
		&Pictorial & Covariance  &\multicolumn{2}{c}{${\rm Cov}_k(X_1,X_2)$} \\\cline{4-5}
		Case &Representation\footnote{The pictorial representation of covariance matrix patterns follows the proposal of \cite{Filzmoser.et.al.Pictorial:2014}.}& Structure  & $k=2,3$ &  $k=4,5$\\
		%%Ex1 ------------------------------------------------
		\hline
		1&%\begin{center}
		%\vspace{0.1cm}
		\includegraphics[scale=0.37]{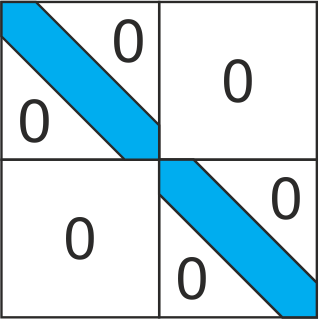}%%%Conf_Diag.png
		%\end{center} 
		&
		\begin{tabular}{l}
			${\rm Cov}(C_1,C_2)=0$\\
			${\rm Cov}(R_1,R_2)=0$\\
		\end{tabular}
		&  $\neq 0$ 	& $=0$\\\hline	%%%\textcolor{purple}{$\times$}%%%\textcolor{purple}{$\checkmark$}
		%%Ex2 --------------------------------------------
		2&
		%\begin{center}
		%\vspace{0.1cm}
		\includegraphics[scale=0.37]{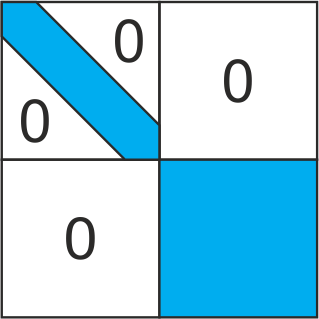}%%%Conf4_SPC.png
		%\end{center} 
		&
		\begin{tabular}{l}
			${\rm Cov}(C_1,C_2)=0$\\
			${\rm Cov}(R_1,R_2)\neq 0$\\
		\end{tabular}
		&  $\neq 0$ 	& $=0$ \\%%%\textcolor{purple}{$\checkmark$}%%%\textcolor{purple}{$\times$}
		\hline	
		%%Ex3 --------------------------------------------
		3&
		 %\begin{center}
		%\vspace{0.1cm}
		\includegraphics[scale=0.37]{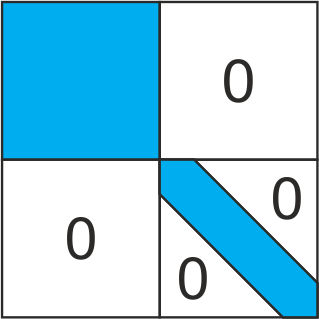} or \includegraphics[scale=0.37]{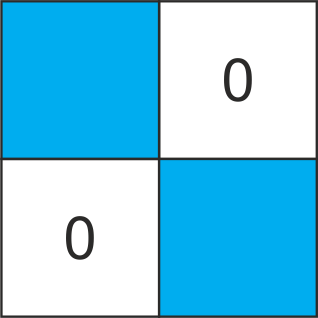}
		%\end{center} 
		&
		\begin{tabular}{l}
			${\rm Cov}(C_1,C_2)=-\delta_k{\rm E}(R_1R_2)$\\%%%{\rm E}(R_1){\rm E}(R_2)
			${\rm Cov}(R_1,R_2)= 0$ or\\ ${\rm Cov}(R_1,R_2)\neq 0$\\
		\end{tabular}
		& $=0$ 	& $\neq 0$ \\%%%\textcolor{purple}{$\times$}%%%\textcolor{purple}{$\checkmark$}
		\hline	
%		%%Ex4 --------------------------------------------
%		 %\begin{center}
%		\vspace{0.1cm}\includegraphics[scale=0.1]{Conf_Diag_Blocks.png}
%		%\end{center} 
%		&
%		\begin{tabular}{l}
%			${\rm Cov}(C_1,C_2)=-\delta_k{\rm E}(R_1R_2)$\\
%			${\rm Cov}(R_1,R_2)\neq 0$\\
%		\end{tabular}
%		& $=0$ \textcolor{purple}{$\times$}	& $\neq 0$ \textcolor{purple}{$\checkmark$}\\
%		
%		\hline
	\end{tabular}
	\label{tab:NullPitfalls}%%%tab:NullSCovPitfalls
\end{table}
\end{minipage}
\vspace{0.2cm}

%%%% CASE 1 %%%%
\begin{example}\label{Ex1:Cov(C1,C2)=0_Cov(R1,R2)=0}
	Let  $X_1$ and $X_2$ be two interval-valued random variables, characterized by $(C_i,R_i)$, $i=1,2$, where $\boma{\mu}$ is the expected value and $\boma{\Psi}$ the covariance matrix of the random vector $(C_1,C_2,R_1,R_2)^t$. %RV 1/11/2020
	In this example, we consider
	\begin{eqnarray}\label{eq:Param_Ex1}
	\boma{\mu}=(0,0,1.5,1.5)^t \;\;\;
	{\rm and}  \;\;\;
	\Psi=\left[\begin{array}{cccc}
	1& \\
	0 & 0.5\\ 
	0 &0 & 0.64\\
	0& 0 & 0    & 0.04\\
	\end{array}\right],
	\end{eqnarray}
	which corresponds to Case 1 of Table \ref{tab:NullPitfalls}, leading to the following symbolic covariance matrices:
	\begin{eqnarray}\label{eq:Param_Ex1cov}
	\boma{\Sigma}_1=\left[\begin{array}{cc}
	1& \\
	0 & 0.5\\ 
	\end{array}\right], 
	&%%%%% Sigma_2
	\boma{\Sigma}_2=\left[\begin{array}{cc}
	1.723& \\
	0.562 & 1.073\\ 
	\end{array}\right], 
	&%%%%% Sigma_3
	\boma{\Sigma}_3=\left[\begin{array}{cc}
	1.241& \\
	0.188 & 0.691\\ 
	\end{array}\right], \nonumber\\ 
	&%%%%% Sigma_4
	\boma{\Sigma}_4=\left[\begin{array}{cc}
	1.723& \\
	0 & 1.073\\
	\end{array}\right], 
	&%%%%% Sigma_5
	\boma{\Sigma}_5=\left[\begin{array}{cc}
	1.241& \\
	0 & 0.691\\ 
	\end{array}\right].
	\end{eqnarray}
	%----------- SCov_k -----------
	%[1,]    1  0.0 NA 1.723 0.562 NA 1.241 0.188 NA 1.723 0.000 NA 1.241 0.000
	%[2,]    0  0.5 NA 0.562 1.073 NA 0.188 0.691 NA 0.000 1.073 NA 0.000 0.691
	%----------- SCor_k -----------
	%[1,]    1    0   NA 1.000 0.414   NA 1.000 0.203   NA     1     0    NA     1     0
	%[2,]    0    1   NA 0.414 1.000   NA 0.203 1.000   NA     0     1    NA     0     1
	
\captionsetup[subfigure]{subrefformat=simple,labelformat=simple,listofformat=subsimple}
\renewcommand\thesubfigure{(\arabic{subfigure})}
\begin{figure}
	%\centering
	\begin{flushleft}
		\captionsetup{justification=raggedright}
		\hspace*{5pt}
		\subfloat[Matrix plots of centers and ranges.]{\label{fig:Ex1_Conv}\includegraphics[width=.5\textwidth]{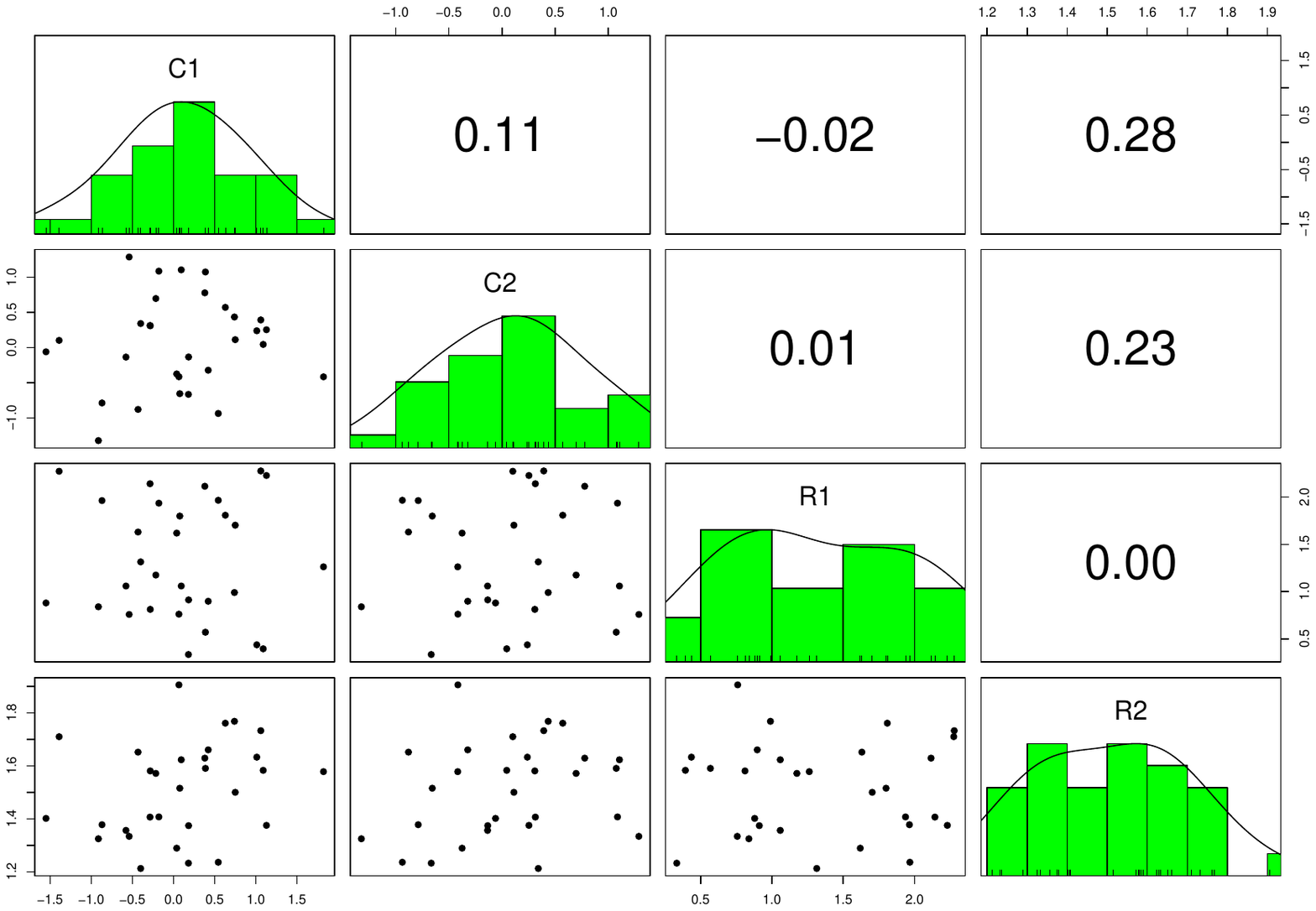}}
		%    \hspace*{5pt}
		\subfloat[Symbolic bivariate scatter plot.]{\label{fig:Ex1_Syn}\includegraphics[width=.5\textwidth]{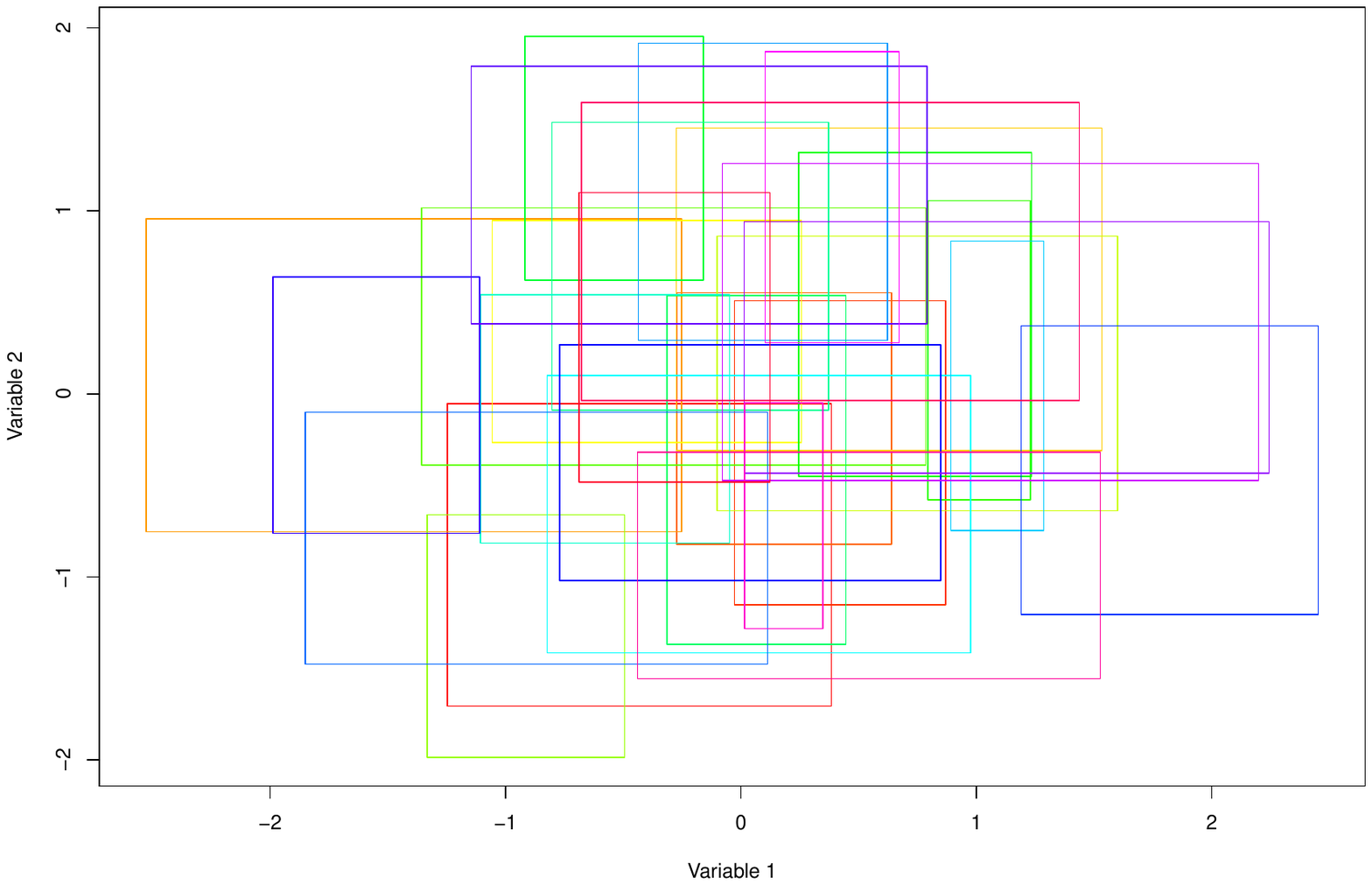}} \    
	\end{flushleft}              
	\caption{Data generated according to Example 1, where the centers and the ranges are uncorrelated.}\label{fig:Null_Scov_Ex1}
\end{figure}

\captionsetup[subfigure]{subrefformat=simple,labelformat=simple,listofformat=subsimple}
\renewcommand\thesubfigure{(\arabic{subfigure})}
\begin{figure}
	%\centering
	\begin{flushleft}
		\captionsetup{justification=raggedright}
		\hspace*{5pt}
		\subfloat[$k=3$.]{\label{fig:Null_Scov_Ex1_A_k=3}\includegraphics[width=.5\textwidth]{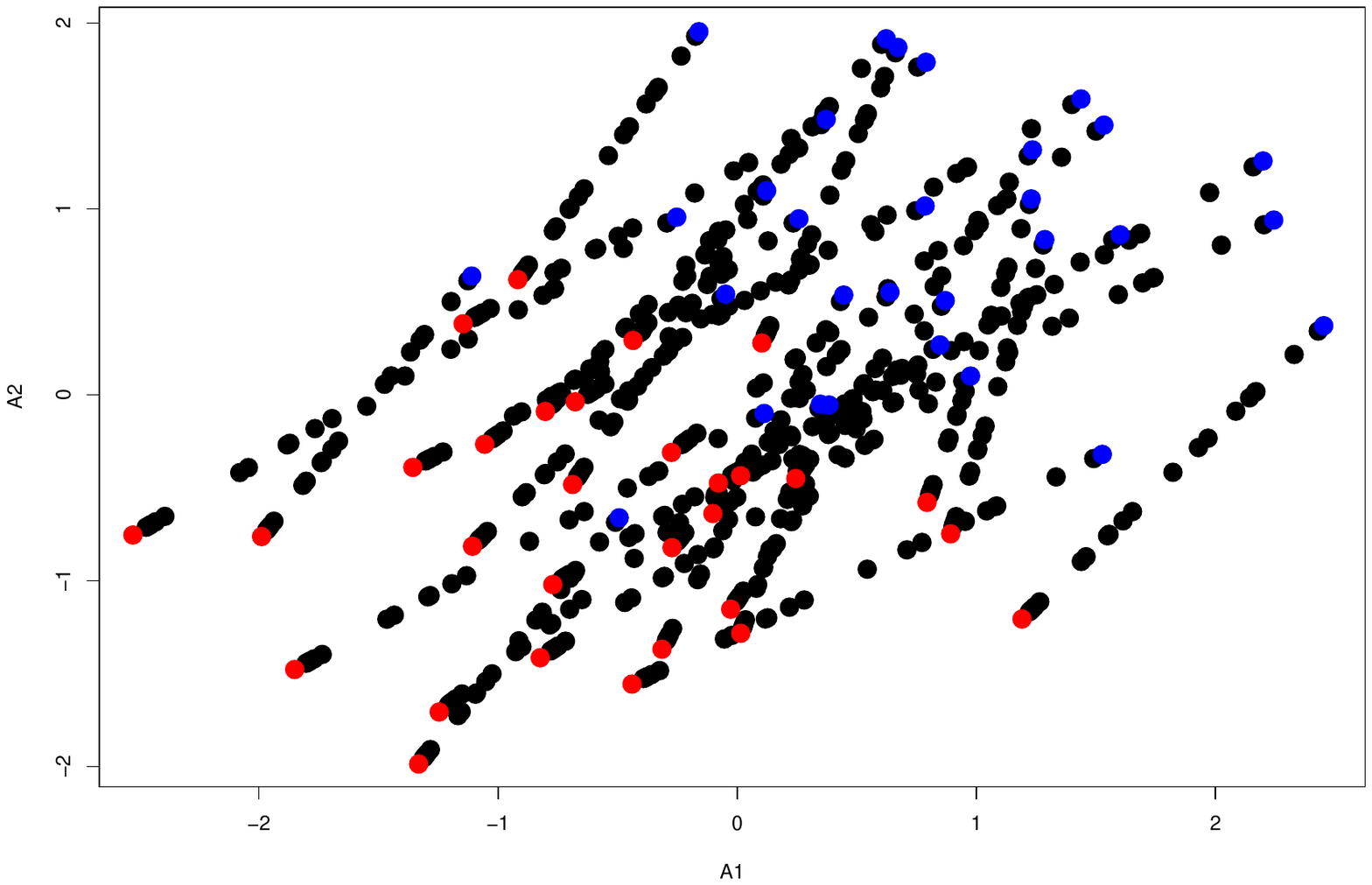}}
		%    \hspace*{5pt}
		\subfloat[$k=5$.]{\label{fig:Null_Scov_Ex1_A_k=5}\includegraphics[width=.5\textwidth]{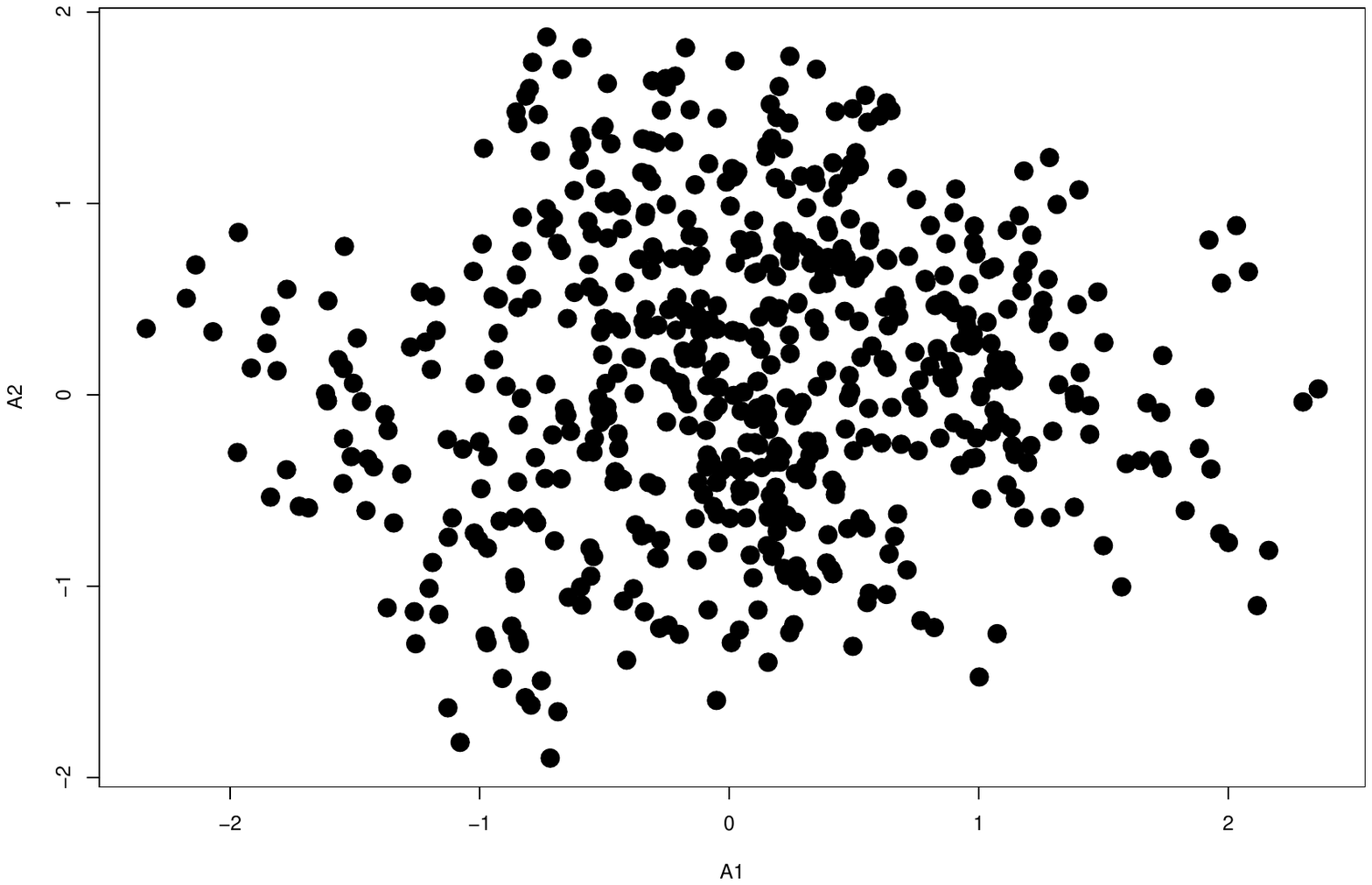}} \    
	\end{flushleft}              
	\caption{Sample of 20 micro-data pairs, $\boma{A}=(A_{1},A_{2})^t$, per symbolic object. Data generated according to Example \ref{Ex1:Cov(C1,C2)=0_Cov(R1,R2)=0}.  In the left plot, the \textcolor{red}{red} (\textcolor{blue}{blue}) dots correspond to the left lower (right upper) vertex of the rectangle that represents each macro-data object.} %%%The centers and ranges follow the distributions specified in Example   \ref{Ex1:Cov(C1,C2)=0_Cov(R1,R2)=0}, which assumes these are non-correlated random variables.%%%Figure \ref{fig:Null_Scov_Ex1_A_k=3} 
	%%%%, where $A_{3j}=C_j+U_3R_j/2$, $A_{5j}=C_j+U_{5j}R_j/2$,  and $U_3$ and $U_{5j}$, are independent variables with uniform distribution over $[-1,1]$, $j=1,2$.
\label{fig:Null_Scov_Ex1_A_k=3,5}
\end{figure}

Clearly, no linear associations between centers and ranges exist, but $\boma{\Sigma}_2$ and $\boma{\Sigma}_3$ show a non-null symbolic correlation (0.414 and 0.203 for definitions 2 and 3, respectively), indicating moderate and weak associations. According to these two definitions and the discussion of Section \ref{sec:Reasoning}, ${\rm Cov}_k(X_1,X_2)={\rm Cov}(A_{1},A_{2})={\rm Cov}(C_1+U\frac{R_1}{2},C_2+U\frac{R_2}{2})$, for $k=2,3$. Thus, since the centers and the ranges are non-correlated random variables, we conclude that the detected association is due to $U$, which is shared by the micro-data pair $(A_{1},A_{2})$. Indeed, the root cause for this problem is the assumption that $U_1 \overset{\rm a.s.}{=}U_2$, which artificially induces a positive association.%RV 2/11/2020

%%%% Simulações de dados e Fig.
To illustrate these findings, we generate a sample of size 30 of centers and ranges from a multivariate normal distribution with the parameters of equation \eqref{eq:Param_Ex1}. The matrix plot and the corresponding symbolic bivariate scatter plot are shown in Figure \ref{fig:Null_Scov_Ex1}. Both plots confirm the absence of linear associations, contradicting the symbolic correlations obtained with definitions 2 and 3 (see equation \eqref{eq:Param_Ex1cov}). To further deepen this issue, we simulated 20 micro-data points according to the reasoning presented in Subsection \ref{subsec:Reasoning_MicroData_Models}, for definitions 3 and 5. The generated  $\boma{A}$ values, represented in Figure \ref{fig:Null_Scov_Ex1_A_k=3,5}, show a strong positive association between the micro-data for $k=3$, and no association for $k=5$.%RV 2/11/2020

\end{example}%%%%%Ex1:Cov(C1,C2)=0_Cov(R1,R2)=0

%----------------------------------------------------
%%% Example 2 - Ex2:Cov(C1,C2)=0_Cov(R1,R2)!=0
%----------------------------------------------------
\begin{example}\label{Ex2:Cov(C1,C2)=0_Cov(R1,R2)!=0}
%%%% CASE 2 %%%%
This example illustrates the Case 2 of Table \ref{tab:NullPitfalls}. Here we consider\\ 
\begin{eqnarray}\label{eq:Param_Ex2}
\boma{\mu}=(0,0,1.3,1.3)^t \;\;\;
{\rm and}  \;\;\;
\Psi=\left[\begin{array}{cccc}
1& \\
0 & 0.5\\ 
0 &0 & 0.16\\
0& 0 & 0.07 & 0.04\\
\end{array}\right],
\end{eqnarray}
a scenario where the centers are uncorrelated but where there is a strong linear association between the two ranges, ${\rm Cor}(R_1,R_2)=0.875$. These choices lead to the following symbolic covariance matrices:%RV 2/11/2020

\begin{eqnarray}%%%%updated
\boma{\Sigma}_1=\left[\begin{array}{cc}
1& \\
0 & 0.5\\ 
\end{array}\right], 
&%%%%% Sigma_2
\boma{\Sigma}_2=\left[\begin{array}{cc}
1.462& \\
0.440 & 0.933\\ 
\end{array}\right], 
&%%%%% Sigma_3
\boma{\Sigma}_3=\left[\begin{array}{cc}
1.154& \\
0.147 & 0.644\\ 
\end{array}\right], \nonumber\\ 
&%%%%% Sigma_4
\boma{\Sigma}_4=\left[\begin{array}{cc}
1.462& \\
0 & 0.933\\ 
\end{array}\right], 
&%%%%% Sigma_5
\boma{\Sigma}_5=\left[\begin{array}{cc}
1.154& \\
0 & 0.644\\ 
\end{array}\right].
\end{eqnarray}
%%%%----------------- SCov_VELHO-------------------
%     X1.1 X2.1    X1.2  X2.2     X1.3  X2.3     X1.4  X2.4     X1.5  X2.5
%[1,]    1  0.0 NA 1.462 0.440 NA 1.154 0.147 NA 1.462 0.000 NA 1.154 0.000
%[2,]    0  0.5 NA 0.440 0.933 NA 0.147 0.644 NA 0.000 0.933 NA 0.000 0.644
%%%%----------------- Scor_updated------------------
%     X1.1 X2.1      X1.2  X2.2       X1.3  X2.3           X1.4  X2.4        X1.5  X2.5
%[1,]    1    0   NA 1.000 0.377   NA 1.00 0.17   NA     1     0    NA     1     0    NA
%[2,]    0    1   NA 0.377 1.000   NA 0.17 1.00   NA     0     1    NA     0     1    NA
%%%%-------Centers and ranges covariance matrix
%     C1 C2  R1   R2
%[1,]  1 0.0 0.00 0.00
%[2,]  0 0.5 0.00 0.00
%[3,]  0 0.0 0.16 0.07
%[4,]  0 0.0 0.07 0.04

\begin{figure}
	%\centering
	\begin{flushleft}
		\captionsetup{justification=raggedright}
		\hspace*{5pt}
		\subfloat[Matrix plots of centers and ranges.]{\label{fig:EX2_Conv}\includegraphics[width=.5\textwidth]{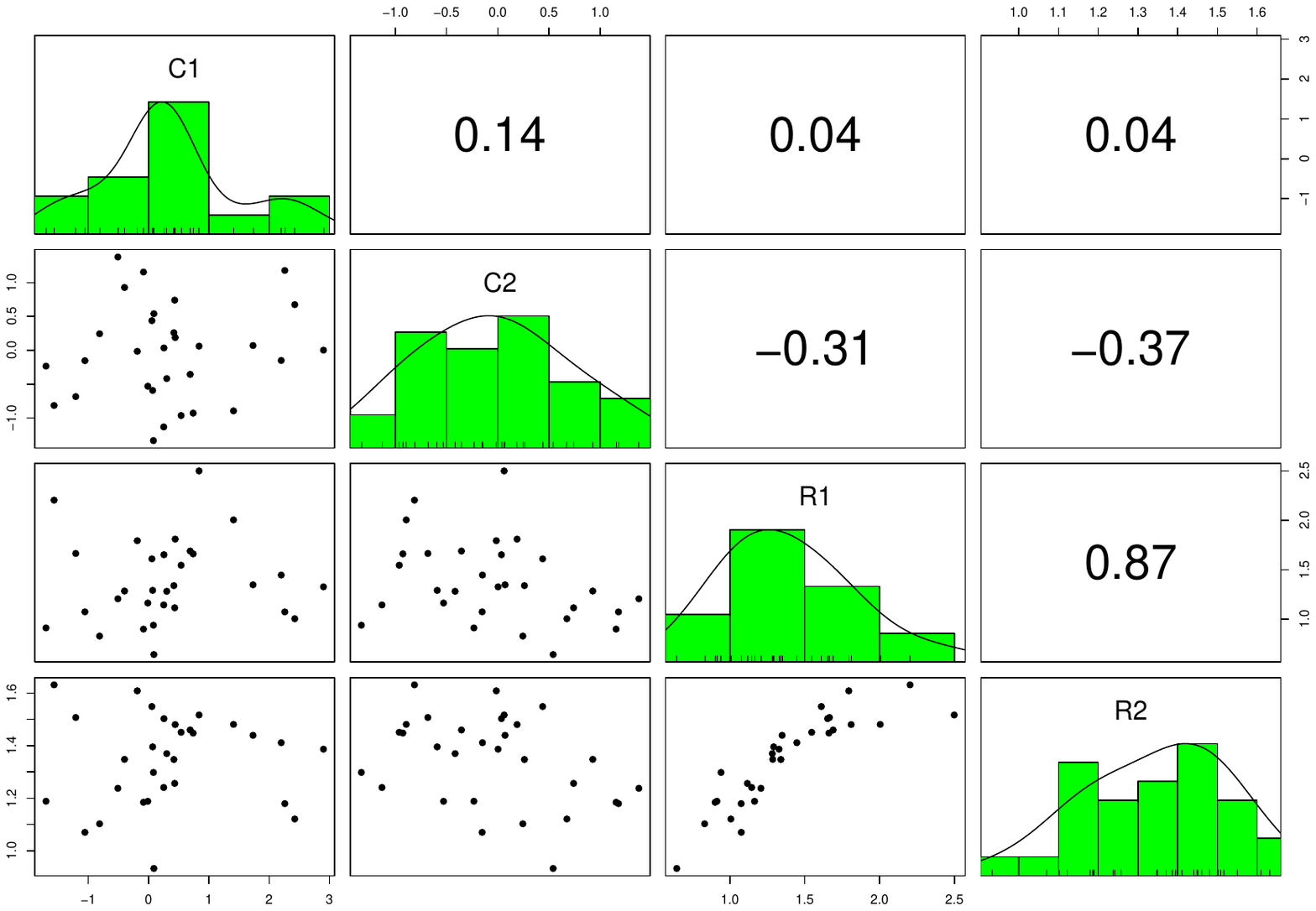}}
		%    \hspace*{5pt}
		\subfloat[Symbolic bivariate scatter plot.]{\label{fig:Ex2_Sym}\includegraphics[width=.5\textwidth]{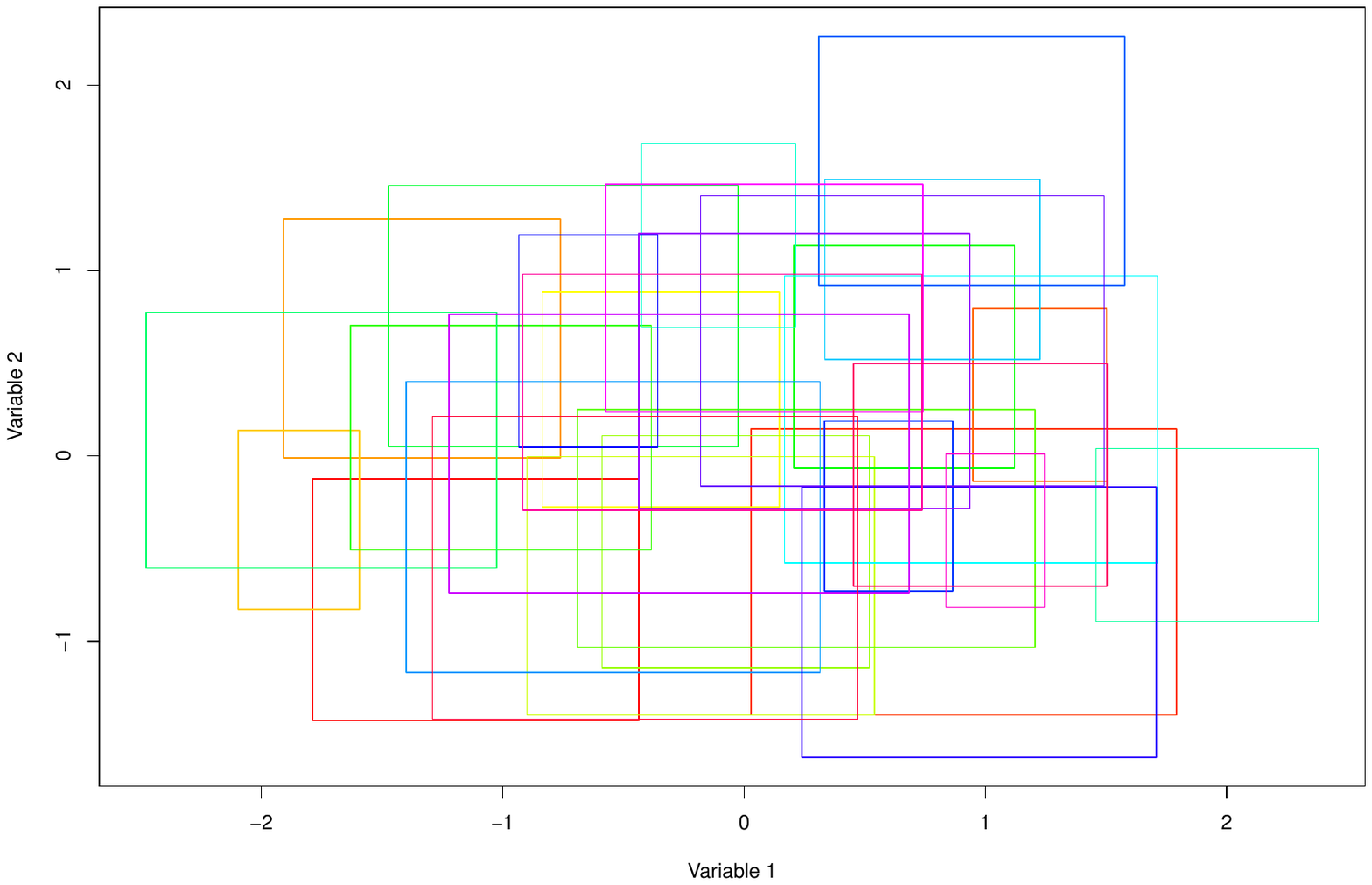}} \    
	\end{flushleft}              
	\caption{Data generated according to Example 2, where centers are uncorrelated and ranges are highly correlated.}
	\label{fig:Null_Scov_Ex2}
\end{figure}

Definitions 1, 4, and 5 all lead to null symbolic correlations, disregarding the strong linear association that exists between the two ranges. Such association is clearly illustrated in Figure \ref{fig:Null_Scov_Ex2}, which corresponds to a sample of size 30 of centers and ranges from a multivariate normal distribution with the parameters of equation \eqnref{eq:Param_Ex2}. The high correlation between the two ranges is apparent in the plot matrix of Figure \ref{fig:EX2_Conv}. Regarding the symbolic bivariate scatter plot of Figure \ref{fig:Ex2_Sym}, the dependence among ranges manifests itself in the relationship between the rectangle lengths. In fact, since we are working with a multivariate normal distribution, then $E(R_2|R_1=r)=0.0840+0.4375r$, meaning that the average ratio between the rectangle lengths should be approximately 0.4375, which is confirmed by Figure \ref{fig:Ex2_Sym}.%RV 2/11/2020

However, even if definitions 2 and 3 lead to non-null symbolic correlation between $X_1$ and $X_2$ (0.377 and 0.170, for $k=2$ and $k=3$, respectively), this is due to the micro-data structure assumed by these definitions, and not to the actual association between the two ranges. In fact, the values would be exactly the same for all eight definitions, would  $\boma{\mu}$ and $\boma{\Psi}$ take the values of equation \eqref{eq:Param_Ex2} except for  $[\Psi]_{3,4}={\rm Cov}(R_1,R_2)=0$. As pointed out before, this is an important limitation of all eight symbolic correlation (covariance) definitions, presented in Table~\ref{tab:reasoningVarCov}.%RV 2/11/2020

\end{example}\label{Ex2:Cov(C1,C2)=0_Cov(R1,R2)!=0}%%%%END Example 2

%----------------------------------------------------
%%% Example 3 - Ex3:Cov(C1,C2)!=0_Cov(R1,R2)=0
%----------------------------------------------------
%%%% CASE 3 %%%%
\begin{example}\label{Ex3:Cov(C1,C2)!=0_Cov(R1,R2)=0}
In this example, we illustrate one of the possible scenarios of  Case 3 of Table \ref{tab:NullPitfalls}. Assuming\\ 
\begin{eqnarray}\label{eq:Param_Ex3}
\boma{\mu}=(0,0,3,3)^t \;\;\;
{\rm and}  \;\;\;
\Psi=\left[\begin{array}{cccc}
1& \\
-9\delta_k & 9\\ 
0 &0 & 0.8\\
0& 0 & 0 & 0.2\\
\end{array}\right],
\end{eqnarray}
leads to the following symbolic covariance matrices:
\begin{eqnarray}
\boma{\Sigma}_1=\left[\begin{array}{rr}
1& \\
0 & 9\\ 
\end{array}\right], 
&%%%%% Sigma_2
\boma{\Sigma}_2=\left[\begin{array}{rr}
3.45& \\
0 & 11.30\\ 
\end{array}\right], 
&%%%%% Sigma_3
\boma{\Sigma}_3=\left[\begin{array}{rr}
1.817& \\
0 & 9.767\\ 
\end{array}\right],\nonumber\\ 
&%%%%% Sigma_4
\boma{\Sigma}_4=\left[\begin{array}{rr}
3.45& \\
-2.25 & 11.30\\ 
\end{array}\right], 
&%%%%% Sigma_5
\boma{\Sigma}_5=\left[\begin{array}{rr}
1.817& \\
-0.752 & 9.767\\ 
\end{array}\right].
\end{eqnarray}
%%%%%----------------------------------------------
%k=1 --> Cov(C1,C2)=0; SCor_1(X1,X2)=0; 
%k=2 --> Cov(C1,C2)=-2.25; Cor(C1,C2)=-0.75; SCor_2(X1,X2)=0; 
%k=3 --> Cov(C1,C2)=-2.25; Cor(C1,C2)=-0.75; SCor_3(X1,X2)=0;
%k=4 --> Cov(C1,C2)= -2.25; Cor(C1,C2)=-0.75; SCor_4(X1,X2)=-0.3603577;
%k=5 --> Cov(C1,C2)= -0.75; Cor(C1,C2)=-0.25; SCor_4(X1,X2)=-0.17923;
%print(round(cbind(SVar1,NA,SVar2,NA,SVar3,NA,SVar4,NA,SVar5),3))
%[1,]    1    0   NA 3.45  0.0   NA 1.817 0.000   NA  3.45 -2.25    NA  1.803 -0.752
%[2,]    0    9   NA 0.00 11.3   NA 0.000 9.767   NA -2.25 11.30    NA -0.752  9.753
%---------------------------------------------------
\begin{figure}
	%\centering
	\begin{flushleft}
		\captionsetup{justification=raggedright}
		\hspace*{5pt}
		\subfloat[Matrix plots of centers and ranges.]{\label{fig:Ex3_Conv}\includegraphics[width=.5\textwidth]{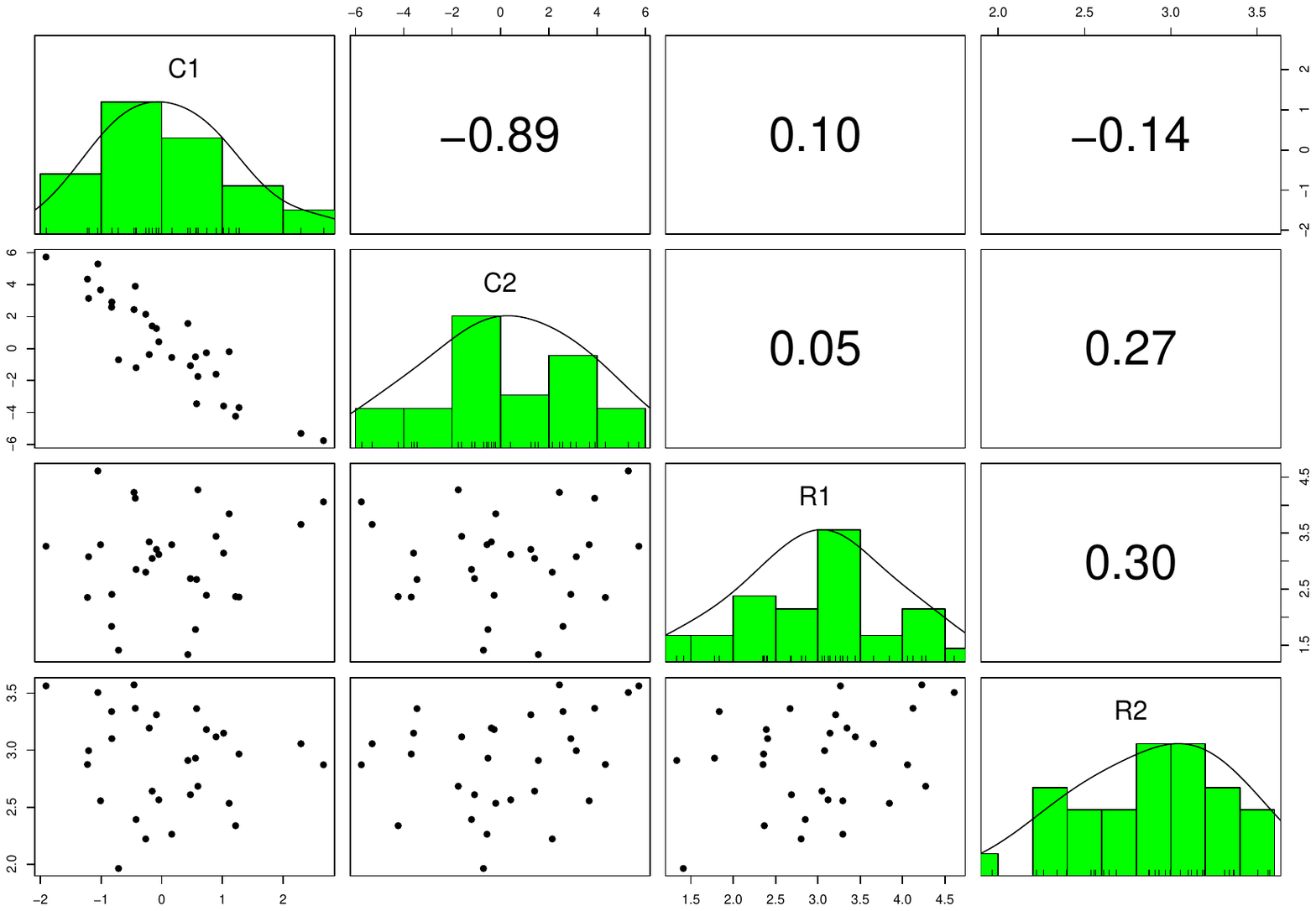}}
		%    \hspace*{5pt}
		\subfloat[Symbolic bivariate scatter plot.]{\label{fig:Ex3_Sym}\includegraphics[width=.5\textwidth]{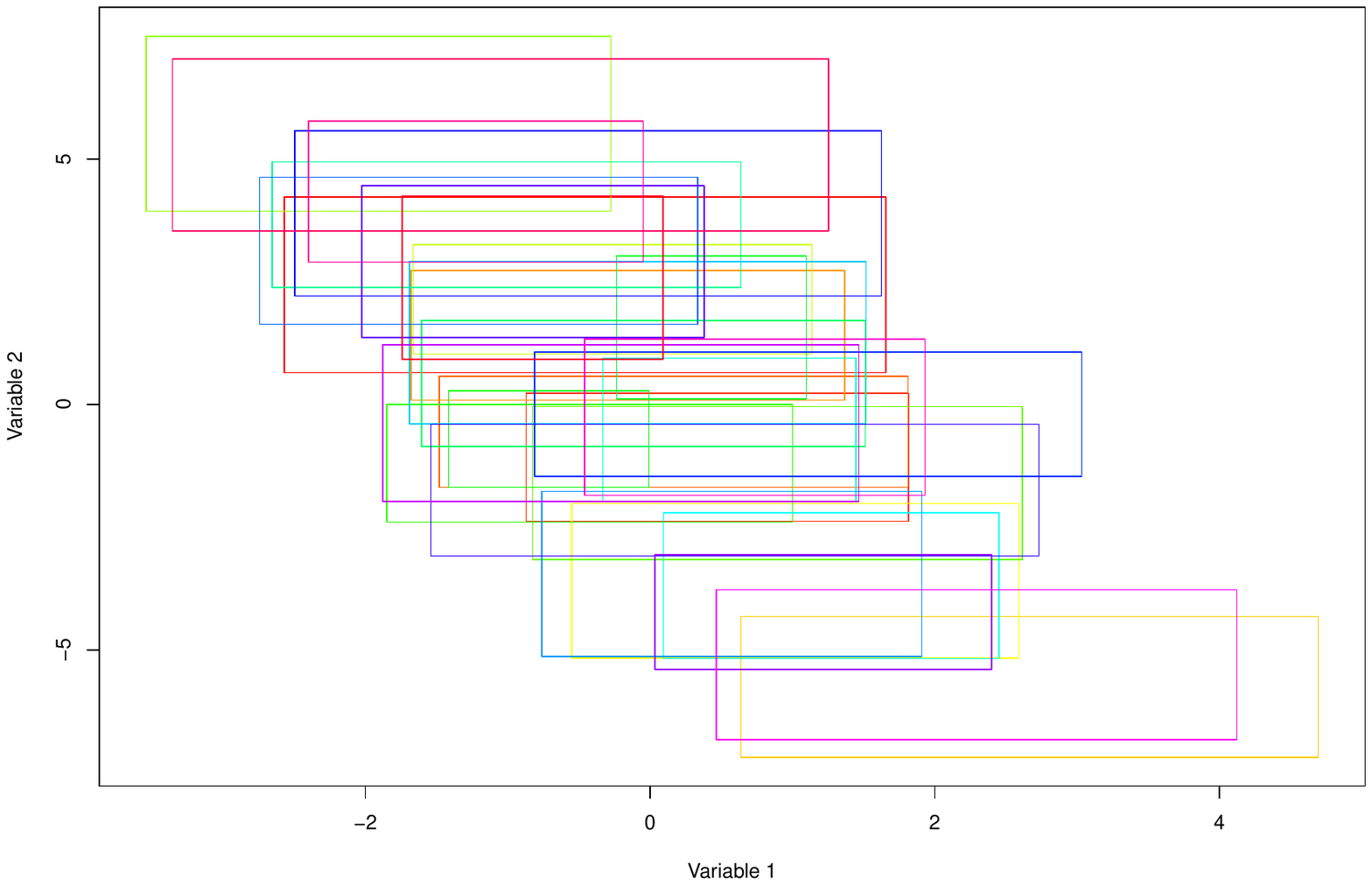}} \    
	\end{flushleft}              
	\caption{Data generated according to Example 3 (for $k=4$), where centers are correlated and ranges are uncorrelated.}
	%%%Centers and ranges are generated from a multivariate normal distribution, with parameters \eqref{eq:Param_Ex3} and $k=4$.
	\label{fig:Null_Scov_Ex3}
\end{figure}

%%%Scenario with correlated centers and ranges, data generated according with definition $k=4$, with null symbolic covariance for $k=2,3$.

The example is built such that there is an association between the two centers which  remains unnoticed by definitions 2 and 3. In these definitions, the contribution of the ranges to the symbolic correlation, which is always positive, can be compensated by a negative correlation between the centers. Specifically, in this example, there is a high negative correlation between the centers (${\rm Cor}(C_1,C_2)=-0.75$, for $k=2$), but this correlation is masked by the positive contribution of the ranges. To confirm this, we generated again 30 observations from a multivariate normal distribution, with the parameters of equation \eqref{eq:Param_Ex3} and $k=4$ ($\delta_k=1/4$). The negative association between the centers is apparent in Figure \ref{fig:Null_Scov_Ex3} but, as mentioned above, this is not captured by definitions 2 and 3.%RV 2/11/2020

\end{example}%%%% CASE 3: Cov(C1,C2)=deltaE(R1)E(R2) and Cov(R1,R2)=0 %%%%

%% file: Rev_Example_SCov_v2.tex
\section{Analysis of datasets}\label{sec:Examples}

In this section, we explore four datasets to illustrate the concepts and results derived previously. The first dataset is the Iris data, which corresponds to four different characteristics of flowers from three different species. The observations are artificially grouped by location, which defines, together with the species, the macro-data object to be studied.%RV 31/10/2020

The next two datasets are Internet traffic data. Internet traffic data is particularly amenable to SDA for two reasons. First, with the increase in Telecommunication's link speeds and in the volume of data stored in Internet sites and repositories (the Big Data problem), it becomes impractical (due to excessive computational costs) to store all measured data (e.g. the arrival time and size of each individual packet observed at a given high-speed link). In this case, the analyst can only have access to summary data (e.g. the maximum and minimum packet sizes observed in a given time period), which can be handled with advantage as a symbolic variables. Second, Internet traffic has several levels/views organized hierarchically (e.g. the packet level, the session level, and the application level) and often, when analyzing one level, there is no interest in having detailed information on the lower levels. For example, when studying TCP connections--session level--we need to have information on the patterns of the packet arrival times and packet sizes of the connection, but not on the arrival time and size of each individual packet that belongs to the connection--packet level.%RV 31/10/2020

The last dataset corresponds to one year monthly credit card expenses of five different types, for which both the micro-data and the macro-data are available.%RV 31/10/2020

\subsection{Iris data}\label{Sec:Iris}

The iris data is probably one of the best known and most used datasets in multivariate analysis, and is also a reference among the SDA community. In the Iris dataset, there is a total of 150 plants, 50 plants from each one of 3 different species (\textit{setosa}, \textit{versicolor}, and \textit{virginica}), and each plant is characterized in terms of \textit{sepal length}, \textit{sepal width}, \textit{petal length}, and \textit{petal width}. Billard and Diday \cite{BillardDiday2006} hypothetically associated sets of five  plants of the same specie to a location. In this case, each macro-data object corresponds to a location of (five) plants from the same specie, and is characterized by four interval-valued variables (sepal length, sepal width, petal length, petal width) ($p=4$). Since there is a total of 150 plants and 5 plants per location, the number macro-data objects is 30 ($n=30$), 10 per specie.%RV 30/10/2020

\begin{figure}
	\centering
	\includegraphics[width=0.8\textwidth]{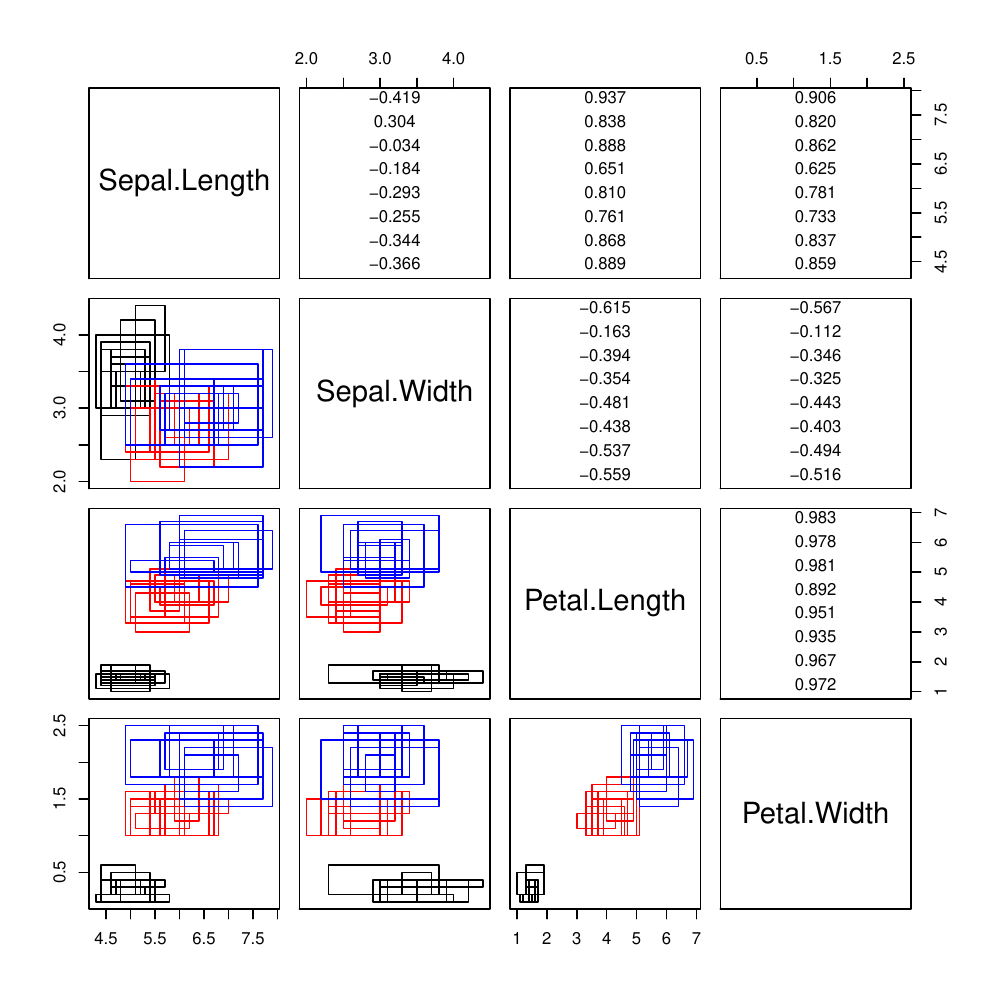}             
	\caption{Symbolic bivariate scatter plots and corresponding sample symbolic correlations of Iris data (the $k$-th row of $(i,j)$-th entry has $\widehat{{\rm Cor}}_k(X_i,X_j)$, $k=1,\ldots,8$). Black rectangles represent locations of setosa, \textcolor{blue}{blue} of versicolor, and \textcolor{red}{red} of virginica.}
	\label{fig:SymbolicPairs_Iris}
\end{figure}

In order to visualize the data, we developed a R \cite{R} function to produce a matrix of all possible combinations of symbolic bivariate scatter plots and the corresponding sample symbolic correlations (available upon request). Figure \ref{fig:SymbolicPairs_Iris} gives an example. The diagonal panels of the matrix show the names of the variables. The lower panels (below the diagonal) show the scatter plots: in the $(i,j)$ panel, $i<j$, each symbolic object is represented by a rectangle whose length and height correspond to the interval values of the $i$-th and $j$-th observations. The upper panels show the sample symbolic correlations: the $(i,j)$ panel, $i>j$, shows the sample correlations that corresponds to the $(j,i)$ scatter plots. Each sample symbolic correlation panel has eight values, each corresponding to one of the definitions introduced in previous sections. Inspired by similar R functions, we call this plot \textit{symbolic pairs}.%RV 30/10/2020

In Figure \ref{fig:SymbolicPairs_Iris}, the locations are colored according to the specie: black for setosa, \textcolor{blue}{blue} for versicolor, and \textcolor{red}{red} for virginica. We notice a strong positive association between petal length and width, and a slightly lower positive association between sepal and petal length, and between sepal length and petal width. In general, locations associated with the setosa specie have smaller center values and inner variability in all variables. Furthermore, versicolor's locations have the highest center values and slightly higher inner variability when compared with virginica's locations, in terms of petal width and length. Similar patterns, although not so clear, are noticed in terms of sepal and petal length. We also notice that petal length and width give a good separation among species. These patterns are also observed in the micro-data.%RV 30/10/2020

\subsection{Internet traffic redirection attacks}

Traffic redirection attacks deviate Internet traffic from its normal routes, allowing attackers interposed between the source and destination of the communications (man-in-middle attackers), to gain access to sensitive information, and/or to degrade the network delay, among other motivations. This type of attack explores vulnerabilities in the BGP protocol, which is responsible for the Internet-wide routing of information. The attack provokes an increase in the round-trip time (RTT) of the communications (from the source to the destination and back to the source), and this feature is usually used to detect the attack. However, the attack is difficult to detect when the attacker is located close to the source or the destination.%RV 31/10/2020

Salvador and Nogueira \cite{Salvador2014} proposed a monitoring infrastructure to detect traffic redirection attacks. This platform comprises a set of geographically dispersed monitoring servers (probes) that periodically measure the RTT to hosts under surveillance (targets). In this case, it may be possible to detect stealthy attacks through suitable statistical processing of the RTTs measured by the various probes.%RV 31/10/2020

To evaluate the RTT, each probe makes 10 RTT measurements every 120 seconds. These RTT measurements form the micro-data, but we have only access to the, average, minimum, and maximum RTT values computed over the 10 individual RTT measurements. The RTT interval (minimum and maximum) forms the macro-data that we study in this example.%RV 28/10/2020

In our example, there are four probes located in Amsterdam, Chicago, Los Angeles (LA), and Johannesburg, and the target is in Hong Kong. During the measurement period, two different attacks were performed, one redirecting the traffic through Moscow and the other through Los Angeles (a site different from the Los Angeles's probe). A total of 2286 RTT interval measurements were collected, where 1799 observations correspond to regular traffic, 241 to traffic redirected through Los Angeles and 246 to traffic redirected through Moscow. For further details on the dataset please see \cite{Subtil2020}.%RV 28/10/2020

\begin{figure}
	\centering
	\includegraphics[width=0.8\textwidth]{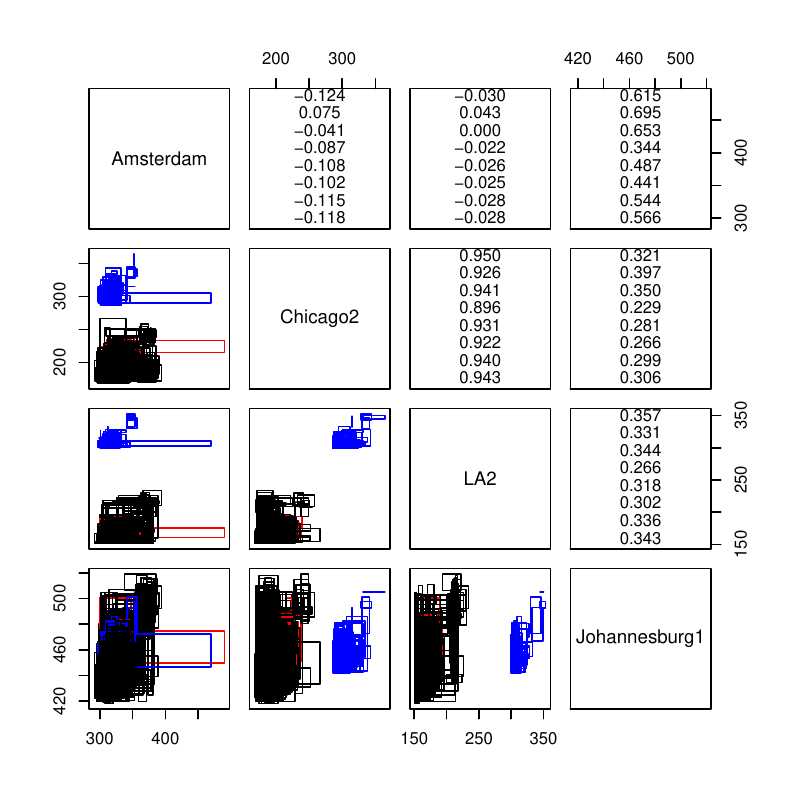}         
	\caption{Symbolic bivariate scatter plots of traffic redirection attacks and respective symbolic correlations. The target is Hong Kong. Black, \textcolor{red}{red}, and \textcolor{blue}{blue}  objects represent regular traffic, traffic redirected through Los Angeles, and traffic redirected through Moscow, respectively.}
	\label{fig:SymbolicPairs_1a8_HongKong}
\end{figure}

By plotting the average RTT versus interval center and by looking to its sample correlation (always higher than 0.979), we found evidence that the micro-data is symmetrically distributed within the intervals, validating one of the assumptions required in the definitions of symbolic covariance matrices. The matrix of scatter plots and associated symbolic correlations based on the eight definitions under study are summarized in Figure \ref{fig:SymbolicPairs_1a8_HongKong}.%RV 31/10/2020 

The bivariate plots seem to indicate that probes located in LA and Chicago separate the regular traffic (black rectangles) from the one redirected through Moscow (\textcolor{blue}{blue}  rectangles). Moreover, this redirection induces a significant increase in RTT values. This justifies the high correlation detected by all sample symbolic definitions of correlation, ranging from 0.896 for definition 4 to 0.950 for definition 1 (see second row and third column of the matrix displayed in Figure \ref{fig:SymbolicPairs_1a8_HongKong}). This Figure also suggests that probes Amsterdam and Johannesburg cannot detect the two attacks. The RTTs measured by these two probes seem to follow a fairly similar pattern, which justifies its medium correlation according to all symbolic definitions.%RV 28/10/2020

The most important observation regarding this example is that, for each of the six probe pairs, the sample symbolic correlation is approximately the same in all eight definitions. This is because all definitions privilege the center of the interval against the range and, in this dataset, the correlation between centers is much stronger than between ranges. Note that in all definitions, the center term has weight one, while the maximum weight associated to the range term is 1/4. To confirm the dominance of centers over ranges, we obtained the sample correlations between the centers:

\begin{small}
	\begin{equation*}
	\widehat{\boma{P}}_{CC}=\left(\begin{array}{rrrrr}
	1.000&\\% -0.124 -0.030  0.615
	-0.124&  1.000\\%  0.950  0.321
	-0.030&  0.950&  1.000\\ % 0.357
	0.615&  0.321&  0.357&  1.000\\ 
	\end{array}\right),
	\end{equation*}
\end{small}
As it can be seen, these correlations are close to the ones of Figure \ref{fig:SymbolicPairs_1a8_HongKong}. For example, the symbolic correlations between Amsterdam and Chicago probes are all close to $[\widehat{\boma{P}}_{CC}]_{2,1}=-0.124$.%RV 31/10/2020

%RV 28/10/2020

\subsection{Internet traffic of backbone networks}\label{sec:IntBackData}

This dataset corresponds to Internet traffic observed in backbone networks (at the core of the Internet), and includes a mixture of Internet applications and attacks, namely Web browsing, file sharing, streaming, video, port scans, and snapshots. It was analyzed in \cite{Oliveira.et.al:2017} using Symbolic Principal Component Analysis for interval-valued data. The dataset comprises 917 traffic objects, corresponding to packets flows of specific applications and attacks, called datastreams. Five traffic characteristics were measured for each datastream, in intervals of 0.1 seconds: (i) number of upstream packets, (ii) number of downstream packets, (iii) number of upstream bytes, (iv) number of downstream bytes, and (v) number of active TCP sessions. The values of these characteristics forms the the micro-data. The macro-data corresponds to the interval (minimum and maximum) values of each traffic characteristic, computed over each datastream.%RV 31/10/2020

As a preprocessing step, we corrected the asymmetry among the micro-data, to allow using the definitions of symbolic covariance matrices presented in section~\ref{sec:Reasoning}.%RV 28/10/2020

In this example, there is a strong association between the centers and ranges of the first four variables. Thus, the definitions that do not take the ranges into consideration when defining the sample covariance lead to low correlations, and the remaining definitions lead to very high correlations. To illustrate this, we consider definitions 3 and 5, where definition 5 does not account for the ranges in the sample covariance. The corresponding sample correlation matrices are:

\begin{small}
	\begin{equation*}
	\widehat{\boma{P}}_{3}=\left(\begin{array}{rrrrr}
 1.000\\% 0.942 0.986 0.893 0.084
 0.942& 1.000\\% 0.917 0.908 0.243
 0.986& 0.917& 1.000\\% 0.915 0.076
 0.893& 0.908& 0.915& 1.000\\% 0.269
 0.084& 0.243& 0.076& 0.269& 1.000\\% 
	\end{array}\right),
	\quad \quad
	\widehat{\boma{P}}_{5}=\left(\begin{array}{rrrrr}
1.000\\% 0.262  0.337 0.205 -0.074
0.262& 1.000\\%  0.239 0.163  0.057
0.337& 0.239&  1.000\\% 0.219 -0.081
0.205& 0.163&  0.219& 1.000\\%  0.068
-0.074& 0.057& -0.081& 0.068&  1.000\\%
	\end{array}\right).
	\end{equation*}
\end{small}
As it can be seen, the sample correlations obtained with definitions 3 and 5 have significantly different values.%RV 29/10/2020 

\subsection{Credit card data}\label{Sec:CreditCard}

This example considers a known symbolic dataset, used in \cite{BillardDiday2003,BillardDiday2006}, that has the merit of having the macro- as well as the micro-data available. The micro-data corresponds to the monthly expenses (in dollars) of three persons, recorded over 12 months, on five different items: food, social entertainment, travel, gas, and clothes. There is a total of $1\,000$ records. The credit card issuer is interested in characterizing the monthly expenses of each person, during one year, over the five different expense types. In our case, each macro-data object corresponds to the expenses of one person in one month, and is characterized by five interval-valued variables, one for each expense type ($p=5$). Since there is a total of 3 persons and 12 months, the number macro-data objects is 36 ($n=36$). The five interval-valued variables are $X_1$, (food), $X_2$ (social entertainment), $X_3$ (travel), $X_4$ (gas), and $X_5$ (clothes).%RV 30/10/2020

\begin{figure}
	\centering
	\includegraphics[width=0.8\textwidth]{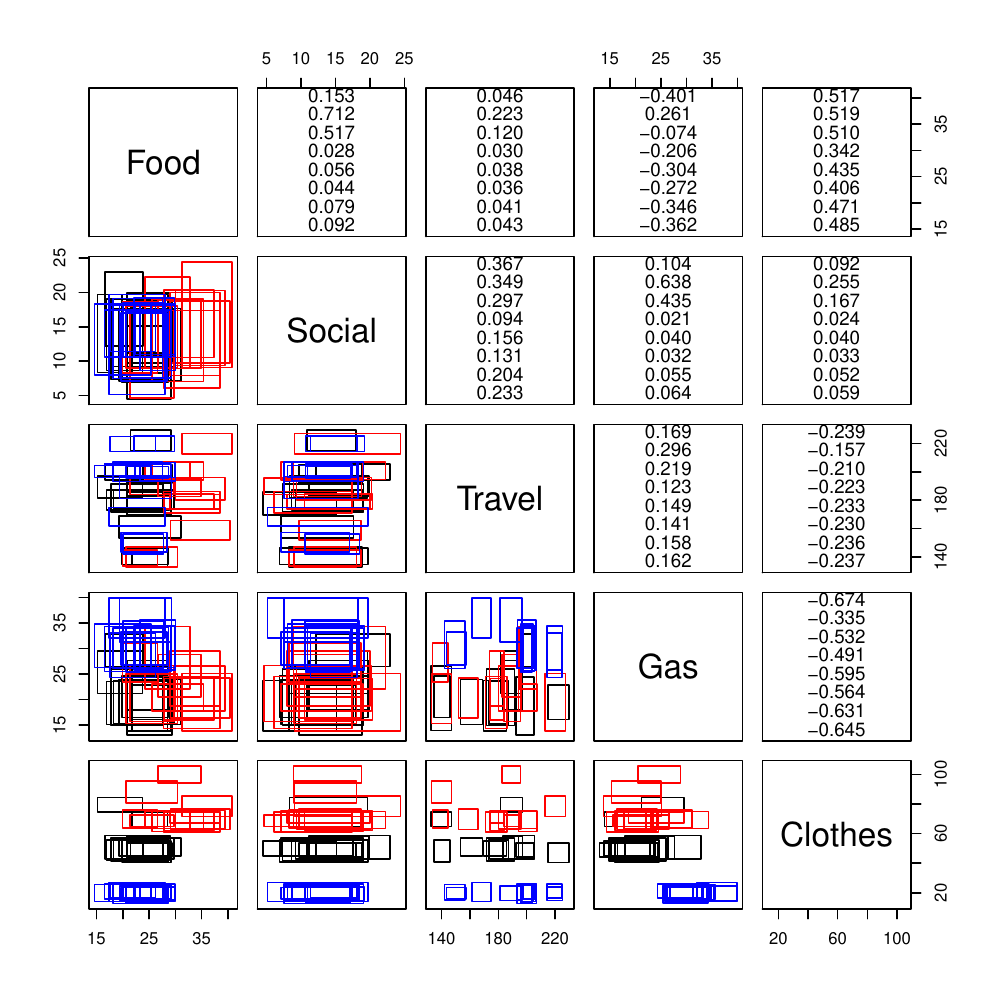}         
	\caption{Symbolic bivariate scatter plots of credit card data and respective symbolic correlations. There are three subjects with monthly expenses measured over a year, colored differently}
	\label{fig:SymbolicPairs_1a8_CreditCards}
\end{figure}

The symbolic bivariate plots and the corresponding sample symbolic correlations are showed in Figure \ref{fig:SymbolicPairs_1a8_CreditCards}, where each person is represented with a different color. It can be said that clothes separates relatively well the three persons. The highest sample symbolic covariances (on absolute value) are the ones between food and clothes (positive values, varying between $\widehat{{\rm Cor}}_4(X_1,X_5)=0.342$ and $\widehat{{\rm Cor}}_2(X_1,X_5)=0.519$) and gas and clothes (negative values, varying between $\widehat{{\rm Cor}}_2(X_4,X_5)=-0.335$ and $\widehat{{\rm Cor}}_1(X_4,X_5)=-0.674$). Thus, subjects with high expenses on clothes tend to have high food and lower gas expenses.%RV 31/10/2020

The results regarding the sample symbolic correlation show that there can be a significant divergence among the eight definitions, which can even have different signs. For example, in the sample correlation between food and gas, seven definitions lead to negative sample correlations (with values between -0.401 and -0.074), and definition $k=2$ leads to a positive sample correlation (0.261). This also happens in the case of Iris data and the sample correlation between sepal length and width, as can be seen in Figure \ref{fig:SymbolicPairs_Iris}.%RV 22/11/2017

In fact, for $k=2$ and $k=3$, the symbolic covariance is a balance between the center's effect, measured by the term ${\rm Cov}(C_j,C_l)$ (which can take positive or negative values), and the range's effect, measured by the term $\delta_k{\rm E}(R_iR_j)$ (which can only take positive values); the weight $\delta_k$ is $1/4$ for $k=2$ and $1/12$ for $k=3$. In this example, $\widehat{{\rm Cov}}(C_1,C_4)=-8.953$ and $\hat{{\rm E}}(R_1R_4)=81.740$. Thus, for $\delta_2=1/4$, the importance of the ranges overcomes the negative covariance between the centers, leading to a positive correlation. Contrarily, for $k=3$, the smaller weight $\delta_3=1/12$ leads to a negative correlation.%RV 30/10/2018

Another clear example of divergence between definitions is the sample symbolic correlation between money spent on food and social entertainment. In this case, the sample correlation ranges from $\widehat{{\rm Cor}}_4(X_1,X_2)=0.028$ to $\widehat{{\rm Cor}}_2(X_1,X_2)=0.712$.%RV 22/11/2017

\begin{figure}
	\centering
	\includegraphics[width=0.6\textwidth]{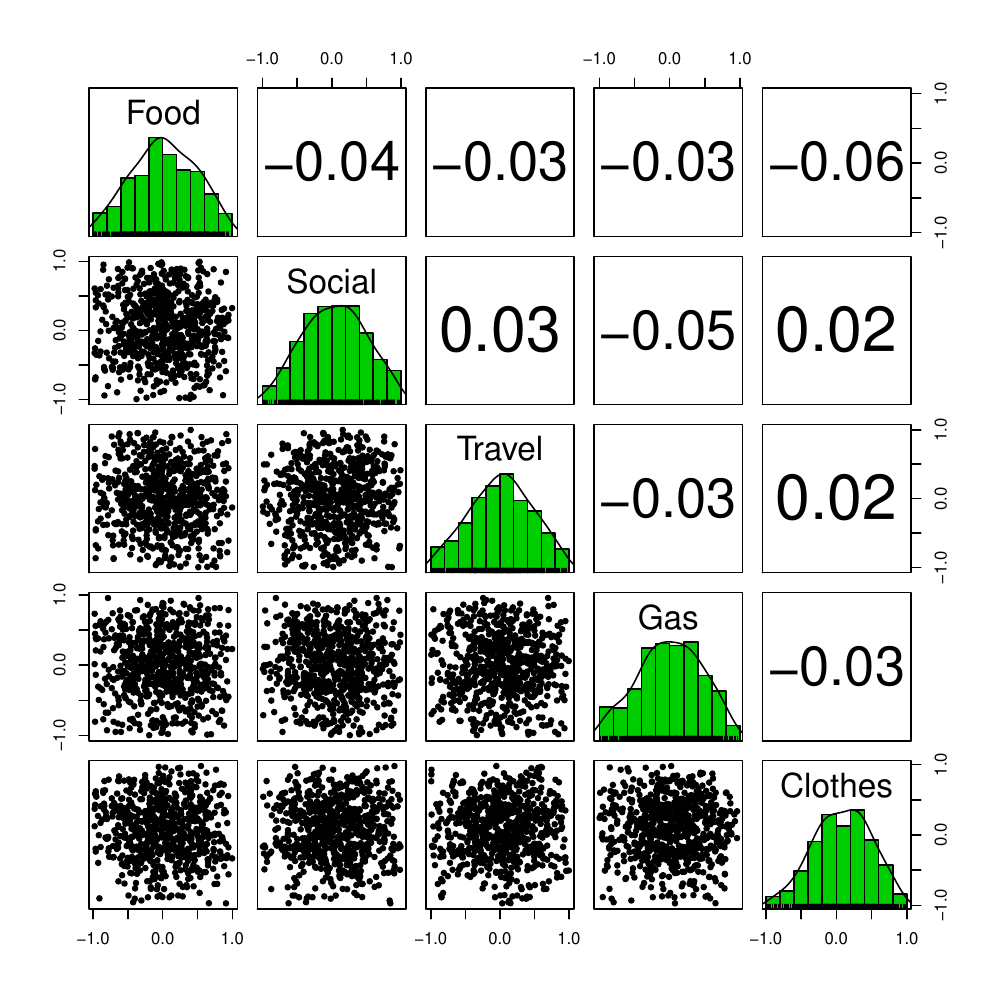}   %            
	\caption{Credit card linearly transformed micro-data, $u_j$, and respective sample correlations. }%%%Microdata on the rectangle borders were eliminated.
	\label{fig:MicroData_CreditCards}
\end{figure}

The divergence between the sample symbolic correlations obtained through the various definitions motivates searching for the most appropriate definition, which can be done since the micro- and macro-data are both available. To address this issue, we represent in Figure \ref{fig:MicroData_CreditCards} the values of the random variables, $U_{j}$, describing the (linearly) transformed micro-data, according to (\ref{Aij}). In this example, the macro-data were not observed directly, but results from the aggregation of the micro-data, according to the credit card issuer criteria. Thus, some of the observed values are used to define the observed interval limits, and to avoid distorting the results, we removed them from the analysis.%Ok 30/10/2020

The scatter plots and sample correlations between concretizations of $U_{j}$ and $U_{l}$ ($j\neq l$) of Figure \ref{fig:MicroData_CreditCards} give indication that these random variables are uncorrelated and that the most promising models among the eight introduced above (see Table~\ref{tab:reasoningVarCov}) are: (i) continuous Uniform, (ii) Triangular, and (iii) truncated Normal distribution, associated with $\boma{\Sigma}_5$, $\boma{\Sigma}_7$, and $\boma{\Sigma}_8$, respectively.%Ok 22/11/2017

We tested these three hypotheses applying goodness of fit tests (we used the Anderson-Darling test \cite{Anderson:2011}) but the null hypothesis was always rejected. This can be explained by the high number of observations, which makes a small departure from the theoretical model statistically significant, even though it might not be significant from the practical point of view.%Ok 22/11/2017

\begin{figure}
	%\centering
	\begin{flushleft}
		\captionsetup{justification=raggedright}
		\hspace*{5pt}
		\subfloat[Food.]{\label{fig:qqPlot_food}\includegraphics[width=.5\textwidth]{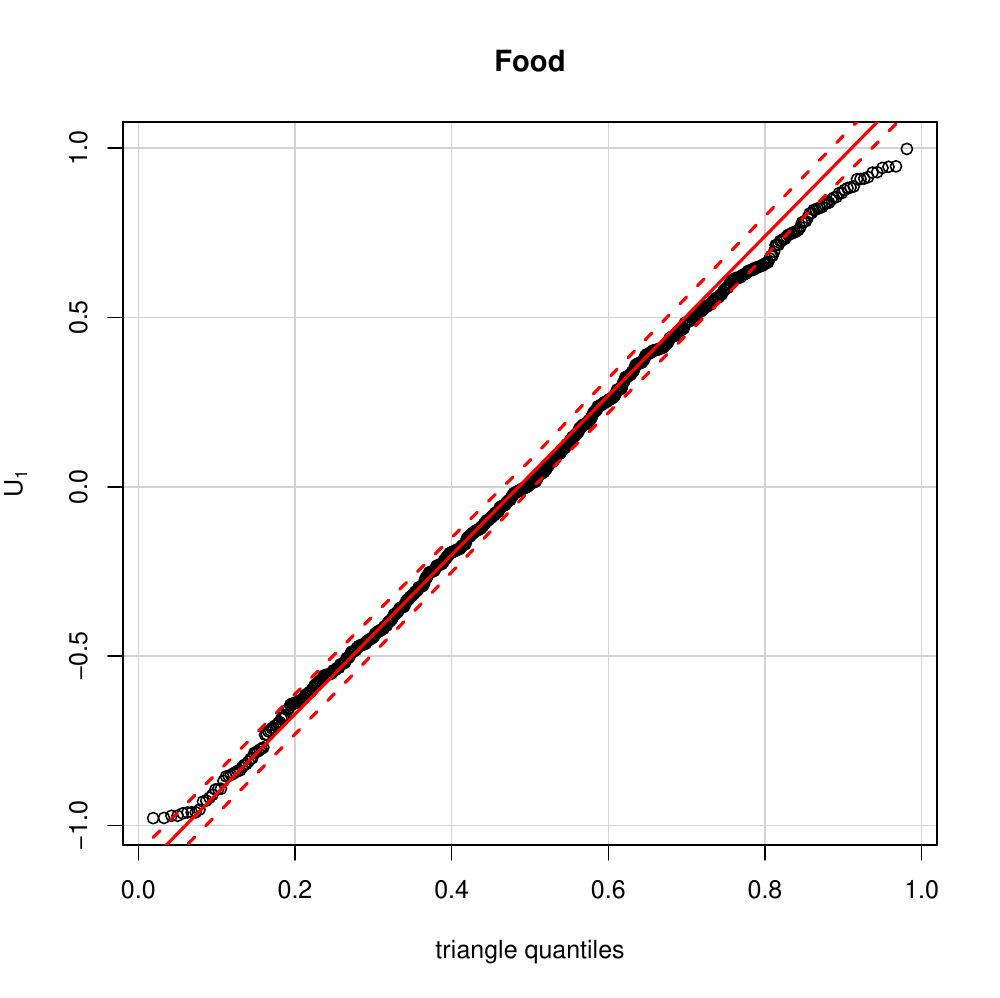}}
		%    \hspace*{5pt}
		\subfloat[Social entertainment.]{\label{fig:qqPlot_social}\includegraphics[width=.5\textwidth]{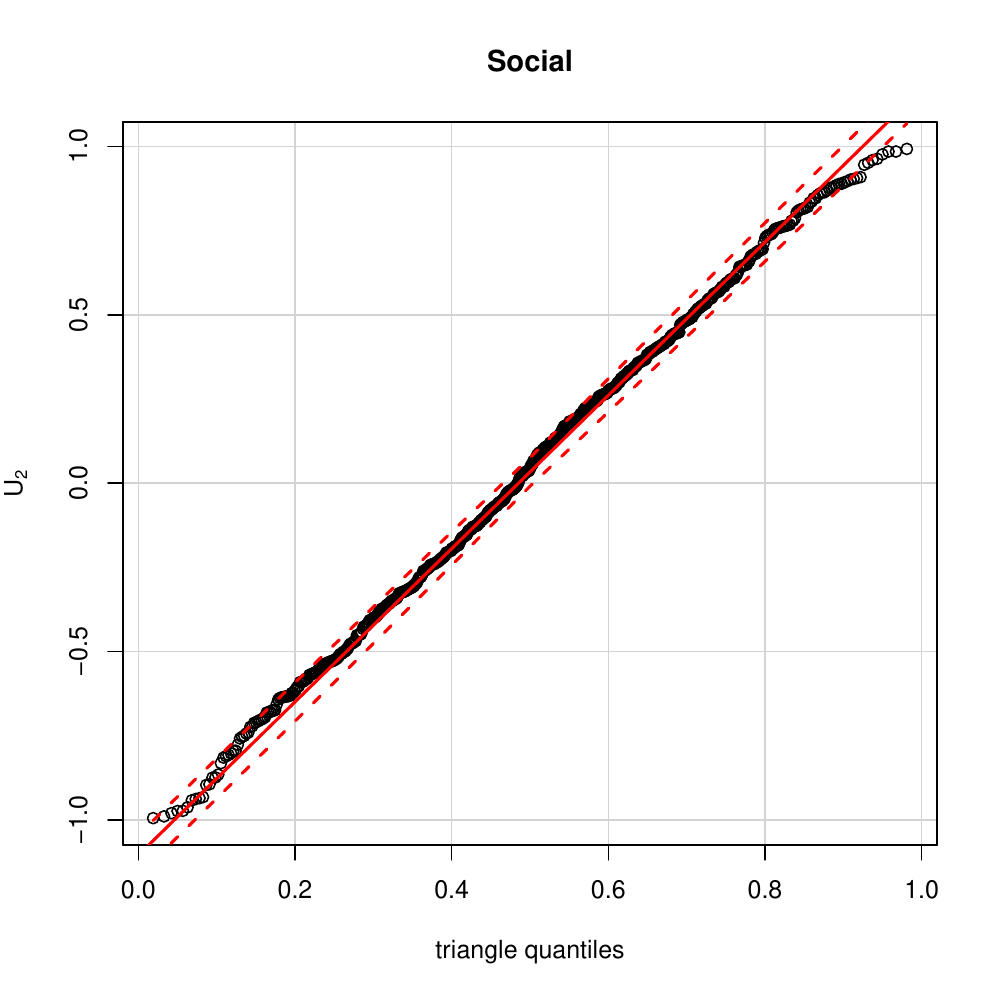}} \    
	\end{flushleft}              
	\caption{Linearly transformed micro-data  quantile-quantile plots, for Triangular[-1,0,1] distribution. Percentage of points outside the 95\% envelope: 9.61\% for  (1) Food ($U_{1}$, worst case) and 1.72\% for (2) Social entertainment ($U_{2}$, best case), where $U_{j}=2(A_{j}-C_j)/R_j$, $j=1,2,\ldots,5$.}
	\label{fig:qqPlot_Food_Social}
\end{figure}

To overcome this problem, we made quantile-quantile plots with 95\% pointwise envelopes (using function \texttt{qqPlot} from package \texttt{car} \cite{car.package:2011} in R), and obtained the percentage of points outside the envelope as a measure of goodness of fit. The Triangular distribution achieved the best results. In this case,
the percentage of points outside the 95\% envelope range between   1.72\% (for social entertainment) to 9.61\% (for food). In Figure \ref{fig:qqPlot_Food_Social}, we show these two quantile-quantile plots, which confirm the goodness of fit. %RV 24/11/2017
From this result, we conclude that definition $k=7$ is the most appropriate for this dataset. Being so, the chosen sample symbolic correlation matrix is
\begin{equation*}
\widehat{\boma{P}}_7=\left(\begin{array}{rrrrr}
1.000 \\%& 0.078 & 0.041 &-0.346 &  0.471\\
0.078 & 1.000 \\%& 0.204 & 0.055 &  0.051\\
0.041 & 0.204 & 1.000 \\%& 0.158 & -0.236\\
-0.346 & 0.055 & 0.158 & 1.000 \\%& -0.631\\
0.471 & 0.051 &-0.236 &-0.631 &  1.000 
\end{array}\right).
\end{equation*}%RV 22/11/2017

For example, there is a medium-sized symbolic positive correlation (0.471) between money spent on food and clothes and a stronger association (even though negative) between gas and clothes (-0.631).%RV 22/11/2017

\captionsetup[subfigure]{subrefformat=simple,labelformat=simple,listofformat=subsimple}
\renewcommand\thesubfigure{(\arabic{subfigure})}

\subsection{Discussion}

The previous examples show that different definitions may lead to quite different results (symbolic sample correlation matrices), depending on the micro-data. Thus, in order to decide on the most appropriate definition, one must have some information about the micro-data structure of the specific case under analysis. Moreover, different definitions may reveal different, but equally interesting, aspects of the data. For example, in \cite{Oliveira.et.al:2017}, which uses the dataset of subsection \ref{sec:IntBackData}, definitions 3 and 5 were used to obtain the sample Symbolic Principal Components and associated scores. In this case, definition 3 highlighted a good separation between the various applications, while definition 5 detected the existence of 3 subgroups in the file sharing application (Torrent) (see Figure 1 of \cite{Oliveira.et.al:2017}).%RV 31/10/2020

We also point out that the existence of different definitions, although it introduces the problem of choosing between them, enriches the set of available statistical tools for exploratory analysis. Moreover, working with assumptions that are difficult or even impossible to validate, is common to many latent variable models. We give two examples. The tetrachoric correlation coefficient \cite{Castellan1966} is a measure of association between two binary random variables assumed to result from the discretization of two latent continuous random variables with normal distribution, and this assumption is impossible to verify. In factor analysis, the common factors are considered random variables with zero mean and unitary variance, and this assumption about the population under study is also impossible to verify.%RV 30/10/2020

%% file: Rev_Conclus_SCov_v3.tex
%\vspace{-25pt}
\section{Conclusions}\label{sec:conclusions}
\vspace{-5pt}

The low cost of information storage combined with recent advances in search and retrieval technologies has made huge amounts of data available, the so-called \textit{big data} explosion. New statistical analysis techniques are now required to deal with the volume and complexity of this data, and Symbolic Data Analysis (SDA) is a promising approach.%RV 2/11/2020

In this paper we propose a model linking the micro-data with the macro-data of interval-valued symbolic variables, which takes a populational perspective. In this model the micro-data is defined as a random vector which is a function of the centers and ranges associated with the macro-data and  random weights characterizing the structure of the micro-data given the associated macro-data. The model defines two scenarios where the various definitions of symbolic covariance matrices already proposed in the literature arise as particular cases. These scenarios correspond to two extreme situations regarding the random weights: in the first scenario the weights are independent random variables and in the second one they are equal variables (almost surely); in both scenarios, the weights are zero mean and uncorrelated latent variables. These conditions on the random weights imply that the current definitions of symbolic covariance matrices rely on micro-data assumptions that may be too stringent, raising applicability concerns. Clearly, more research is required in this area.%RV 2/11/2020

We discuss in detail several cases where the existence/absence of correlations in the macro-data are not correctly captured by the definitions. These inconsistencies are explained by the (too restrictive) underlying micro-data assumptions. These cases also highlights that, in the context of current definitions, a null symbolic covariance can not be interpreted as absence of association. %RV 2/11/2020
This reinforces the need of further research on how to measure associations between interval-valued variables. %%RO 2Nov2020

The analysis of four different datasets further explores the various definitions of symbolic covariance matrices. We show that, when using real data, there can be a large divergence between the various definitions, in particular when there is a strong association between the ranges in the data. Thus, in order to select the most appropriate definition, one must have some knowledge about the micro-data structure. For datasets where both the micro- and macro-data are available we were able to select the definition that better explains the data. We also highlight that different definitions may reveal different aspects of the data, which can be used in exploratory data analysis.%RV 2/11/2020